\documentclass[notitlepage, 11pt]{article}

\usepackage{color}
\usepackage{latexsym}
\usepackage{amsmath, amssymb, amsfonts}

\newtheorem{lemma}{Lemma}
\newtheorem{theorem}{Theorem}

\newtheorem{definition}{Definition}

\newtheorem{proposition}{Proposition}

\newtheorem{assumption}{Assumption}

\DeclareMathOperator*{\argmin}{arg\,min}

\def\thelemma{\arabic{section}.\arabic{lemma}}
\def\thetheorem{\arabic{section}.\arabic{theorem}}
\def\thecorollary{\arabic{section}.\arabic{corollary}}
\def\thedefinition{\arabic{section}.\arabic{definition}}
\def\theexample{\arabic{section}.\arabic{example}}
\def\theproposition{\arabic{section}.\arabic{proposition}}

\def\theassumption{\arabic{section}.\arabic{assumption}}

\def\theremark{\arabic{section}.\arabic{remark}}

\newcommand{\manualnames}[1]{
\def\thelemma{#1.\arabic{lemma}}
\def\thetheorem{#1.\arabic{theorem}}
\def\thecorollary{#1.\arabic{corollary}}
\def\thedefinition{#1.\arabic{definition}}
\def\theexample{#1.\arabic{example}}
\def\theproposition{#1.\arabic{proposition}}
\def\theassumption{#1.\arabic{assumption}}
\def\theremark{#1.\arabic{remark}}
}

\newcommand{\beginsec}{
\setcounter{lemma}{0}
\setcounter{theorem}{0}
\setcounter{corollary}{0}
\setcounter{definition}{0}
\setcounter{example}{0}
\setcounter{proposition}{0}
\setcounter{condition}{0}
\setcounter{assumption}{0}
\setcounter{conjecture}{0}
\setcounter{problem}{0}
\setcounter{remark}{0}
}

\newcommand{\noi}{\noindent}

\newcommand{\Uf}{\mathfrak{U}}
\newcommand{\E}{\mathbb{E}}
\newcommand{\R}{\mathbb{R}}

\newcommand{\p}{\mathbb{P}}
\newcommand{\N}{\mathbb{N}}
\newcommand{\I}{\mathcal{I}}
\newcommand{\Ir}{\mathbb{I}}
\newcommand{\Jr}{\mathbb{J}}
\newcommand{\Kr}{\mathbb{K}}
\newcommand{\Lr}{\mathbb{L}}

\newcommand{\la}{\lambda}
\newcommand{\sig}{\sigma}

\newcommand{\eps}{\varepsilon}

\def\limn{\underset{n\rightarrow\infty}{\lim}} 
\newcommand{\Barr}{{\beta_0}}

\newcommand{\ee}{C}
\newcommand{\cc}{s}
\newcommand{\cco}{s_1}
\newcommand{\cct}{s_2}

\newcommand{\T}{\bar T}

\newcommand{\ph}{\varphi}
\newcommand{\vr}{\varrho}
\newcommand{\al}{\alpha}
\newcommand{\gam}{\gamma}

\newcommand{\s}{\sigma}
\newcommand{\del}{\delta}
\newcommand{\om}{\omega}

\newcommand{\Gam}{\mathnormal{\Gamma}}
\newcommand{\Del}{\mathnormal{\Delta}}

\newcommand{\La}{\mathnormal{\Lambda}}
\newcommand{\X}{\mathnormal{\Xi}}

\newcommand{\Om}{\mathnormal{\Omega}}

\newcommand{\Z}{{\mathbb Z}}

\newcommand{\EE}{{\mathbb E}}

\newcommand{\PP}{{\mathbb P}}

\newcommand{\brho}{\boldsymbol{\rho}}
\newcommand{\by}{\boldsymbol{y}}
\newcommand{\bB}{\boldsymbol{B}}
\newcommand{\bS}{\boldsymbol{S}}
\newcommand{\bT}{\boldsymbol{T}}

\newcommand{\calC}{{\cal C}}
\newcommand{\calD}{{\cal D}}
\newcommand{\calE}{{\cal E}}
\newcommand{\calF}{{\cal F}}
\newcommand{\calG}{{\cal G}}
\newcommand{\calH}{{\cal H}}
\newcommand{\calI}{{\cal I}}

\newcommand{\calL}{{\cal L}}

\newcommand{\calP}{{\cal P}}

\newcommand{\calX}{{\cal X}}

\newcommand{\AC}{{\cal AC}}

\newcommand{\oo}{\overline}
\newcommand{\skp}{\vspace{\baselineskip}}

\newcommand{\w}{\wedge}

\newcommand{\pl}{\partial}

\newcommand{\dist}{{\rm dist}}

\newcommand{\iy}{\infty}
\newcommand{\osc}{{\rm osc}}
\newcommand{\ds}{\displaystyle}
\newcommand{\A}{{\cal A}}
\newcommand{\B}{{\cal B}}
\newcommand{\IA}{{\it IA}}
\newcommand{\ST}{{\it ST}}
\newcommand{\qed}{\hfill $\Box$}

\newcommand{\trian}{\mathnormal{\Delta}}

\newcommand{\uzero}{0}
\newcommand{\uh}{{\bar h}}
\newcommand{\ur}{{\bar r}}

\newcommand{\uD}{{\bar D}}
\newcommand{\ue}{e}

\newcommand{\upsi}{{\bar\psi}}
\newcommand{\uxi}{{\bar \xi}}
\newcommand{\utheta}{{\theta}}

\newcommand{\hs}{0}

\title{Asymptotically optimal control for a multiclass queueing model in the moderate deviation
heavy traffic regime\thanks{Research supported in part by the ISF (Grant 1315/12)}}
\author{Rami Atar\thanks{Department of Electrical Engineering,
Technion--Israel Institute of Technology,
Haifa 32000, Israel}
\and
Asaf Cohen\thanks{Department of Mathematics,
University of Michigan,
Ann Arbor, MI 48109, USA. Email: shloshim@gmail.com, web: https://sites.google.com/site/asafcohentau/}
}

\date{\today}

\begin{document}
\maketitle

\begin{abstract}
A multi-class single-server queueing model with finite buffers,
in which scheduling and admission of customers are subject to control, is studied
in the moderate deviation heavy traffic regime. A risk-sensitive cost set
over a finite time horizon $[0,T]$ is considered.
The main result is the asymptotic optimality of a control policy
derived via an underlying differential game.
The result is the first to address a queueing control problem at
the moderate deviation regime that goes beyond models having the so called
pathwise minimality property.
Moreover, despite the well known fact that an optimal control over a finite time interval
is generically of a nonstationary feedback type,
the proposed policy forms a stationary feedback, provided $T$ is
sufficiently large.

\skp

\noi

\skp

\noi

\end{abstract}

\section{Introduction}\label{sec1}
\beginsec

This paper continues a line of research started
in \cite{ata-bis} that aims at analyzing queueing control problems (QCPs) at
the moderate deviation (MD) heavy traffic regime.
The model under consideration consists of a server that serves
customers from a number of classes, where allocation
of the effort among classes is dynamically controlled.
Customers are kept in buffers of
finite size, one buffer for each class, and those that arrive to find a full buffer are lost.
It is also possible to reject arrivals when buffers are not full.
This control system is considered with a risk-sensitive (RS) cost,
that accounts for holding of customers in the buffers as well as for rejections.
At the heart of the analysis lies a differential game (DG) that has been analyzed
in \cite{ata-coh}. This paper proves the validity of the prediction of \cite{ata-coh} that
the DG governs the scaling limit, by showing that the QCP's value converges
to the DG's value, and identifying an asymptotically optimal (AO) policy
for the former that is constructed based on the latter.
The limit result for the model treated in \cite{ata-bis} was built
on {\it pathwise minimality}, a property that considerably simplifies
the analysis, which does not hold in our setting. Instead, the proof here is based
on the Bellman (or the dynamic programming) equation and, specifically,
a {\it free boundary point} characterized by it governs the asymptotic behavior.

Traditionally, heavy traffic analysis of queueing models, and particularly QCPs, is carried
out under the regime of {\it diffusion-scale deviation} (sometimes also
referred to as {\it ordinary deviation}), but it
is also relevant at the MD scale, where relatively
few results exist \cite{Puhal-1999}, \cite{Majewski06}, \cite{ganesh}.
The roots of {\it large deviation} (LD) analysis of control systems go back to
Fleming \cite{fle-1}, who studies the associated Hamilton-Jacobi equations.
The connection of RS cost to DG was made by Jacobson \cite{Jac-1}.
Analyzing RS control by LD tools and the formulation of the corresponding
maximum principle are due to Whittle \cite{whittle-1}.
Various aspects of this approach have been studied
for controlled stochastic differential equations,
for example in \cite{dup-mce}, \cite{fle-mce}, \cite{fle-sou}.
The treatment of a QCP at the MD scale is similar to that at the LD scale
(in papers such as \cite{ata-dup-shw}) as far as the tools are concerned,
but there are reasons to believe that the games obtained in the MD regime are solvable
more often than in the LD regime.
This statement is supported by the fact that the paper \cite{ata-bis} solves
a DG for the MD scale, whereas
a solution of an analogous DG for the LD regime is not known in general
(see \cite{ata-men} for a partial solution of the latter,
and an open problem regarding its general solution).
Similarly, the DG of this paper has been solved in \cite{ata-coh}, but an explicit
solution for the LD analogue is not known.
An additional advantage the MD regime has over LD is
the invariance to the stochastic data, specifically,
the arrival and service time distributions, as long as they possess
exponential moments, where LD results are more sensitive.
The combination of these properties provides a great deal of motivation for working
at the MD scale.

The aforementioned pathwise minimality property has been the basis for solving QCPs
in diffusion scale heavy traffic asymptotics in various works in the past
(for example, \cite{har-lop}).
To describe this property, consider the simple diffusion control problem of minimizing
a cost $J(\zeta)$
over all control processes $\zeta$ having $\R_+$-valued nondecreasing sample paths.
The cost takes the form
\[
J(\zeta)=\E\int_0^Th(\xi_t)dt,
\]
where $h:\R_+\to\R_+$ is nondecreasing,
$\xi_t=x+w_t+\zeta_t$, $x\ge0$ and $T>0$ are fixed,
$w$ is a standard Brownian motion, and the constraint $\xi_t\ge0$ for all $t\ge0$
must be met. The solution is to set $\zeta_t=-\inf_{s\le t}[(x+w_s)\w0]$, making
$\xi$ a reflected Brownian motion starting from $x$, reflecting at $0$.
This follows by the well-known
fact that a.s., for all $t$, $x+w_t+\zeta_t\le x+w_t+\tilde\zeta_t$, for any control
$\tilde\zeta$ meeting the constraint (see, for example, Section 2 of \cite{cheman}).
Although this problem is simpler than typical diffusion control problems
in the literature, pathwise solutions of these problems owe to this simple
property (or sometimes multidimensional versions thereof).

The DG of \cite{ata-bis}, identifying the MD asymptotics of a QCP,
was also solved by such a consideration.
A simplified version of this game, presented with
one-dimensional instead of multidimensional dynamics, is as follows.
It is a zero-sum game with payoff
\begin{equation}\label{104}
\tilde J(\tilde\la,\tilde\mu,\zeta)=\int_0^Th(\ph_t)dt-\int_0^T(a\tilde\la_t^2 +b\tilde\mu_t^2)dt,
\end{equation}
and dynamics
\begin{equation}\label{105}
\ph_t=x+\int_0^t(\tilde\la_s-\tilde\mu_s)ds+\zeta_t.
\end{equation}
Here, $x$, $T$, $a$ and $b$ are positive constants, and $h$ is again nondecreasing.
The function $\zeta$ is a control for the minimizing player, taking values
in $\R_+$ and is nondecreasing, while
the functions $\tilde\la$ and $\tilde\mu$ form a control
for the maximizing player, and are nonnegative. The constraint $\ph_t\ge0$ for all $t$
must be satisfied.
In this game, the functions $\ph$ and $\zeta$ represent MD-scaled queue length and
idleness processes, while $\tilde\la$ and $\tilde\mu$ stand for MD-scaled perturbations
of the arrival and service processes. The function $h$ is the running cost in the underlying
RS cost, while the second term in \eqref{104}
corresponds to penalty associated with changes of measure, and its form originates
from the LD rate function (background on the structure of DGs governing
RS control asymptotics appears in \cite{Fleming2006b} and \cite{Fleming2006}).

It is easy to see how pathwise minimality can be used
once again to find an optimal strategy for the minimizer.
Namely, for $(\tilde\la,\tilde\mu)$ given, setting $\zeta_t=-\inf_{s\le t}(\psi_s\w0)$,
where $\psi_t=x+\int_0^t(\tilde\la_s-\tilde\mu_s)ds$, results with
$\ph$ that bounds from below any other dynamics adhering
to the constraint.
Significantly, this pathwise minimality property provides not only a solution to
the game but also the basis of the AO proof
in \cite{ata-bis}, as one can mimic the behavior of this strategy
to come up with a policy for the queueing model that is automatically AO.

It turns out that one cannot argue along the same lines for the game
obtained under LD scaling. Indeed, note carefully that
the solution method just presented
uses the fact that the second term in \eqref{104} does not
depend on the control for the minimizing player. However,
under LD scaling, the corresponding penalty term, accounting for changes of measure,
involves controls of both players. This makes it impossible to obtain
a pathwise solution in the same fashion.
This point is explained in detail in Section 1 of \cite{ata-bis}.

Although pathwise minimality is useful when it applies,
it is not generic even under the diffusion and MD regimes.
A natural approach to handle more general settings,
that has been used in numerous papers on diffusion scale asymptotics,
is to appeal to dynamic programming methods to solve diffusion control problems and
then use these solutions as a vehicle for analyzing
the QCP (for a small sample of these papers see \cite{AMR}, \cite{BG2012}, \cite{war-kum}).
This approach has not been considered before for MD asymptotics of QCPs.
The model studied in this paper is indeed suitable for such an approach, and in fact
constitutes a prototype for QCPs that are too complex to possess directly solvable DGs,
while a solution via dynamic programming is available.

The DG for our model differs from the one presented
above. Again, we present it in a slightly simplified way; the precise details
appear in Section \ref{sec2}.
The payoff \eqref{104} has an additional term $\vr_T$,
the dynamics \eqref{105} has an additional term $-\vr_t$, and the constraint
$\ph_t\ge0$ is strengthened to $\ph_t\in[0,D]$ for all $t$, where $D$ is a constant.
The $\R_+$-valued nondecreasing function $\vr$ represents cumulative rejections,
and is considered part of the control
for the minimizing player; that is, in this game the minimizing player
controls the pair $(\zeta,\vr)$.
The constraint stems from the finiteness of the buffer,
and the constant $D$ is related to the buffer size (in fact,
it is the buffer size measured in units of MD-scaled workload).
In \cite{ata-coh}, this game was analyzed via a free boundary value problem,
and solved for the value function and optimal strategy.
The contribution of this paper is to substantiate the relation
of the queueing model to the game in a rigorous manner, showing that the latter
indeed governs the MD asymptotics.
This is established by proving that the value of the RS QCP
converges to that of the DG, and translating the DG's optimal strategy
into an AO policy for the QCP.

The proof of convergence of the RS value to the DG value is performed
in two main steps, namely bounding the former from below and from above by the latter.
We refer to them as the lower and upper bound, respectively.

The proof of an asymptotic lower bound of a RS
cost in models that do not have control is often based on the classical
proof of (the lower bound in) Varadhan's lemma \cite{DZ} for a sequence
of processes satisfying the large deviations principle.
This is the case for example for queueing models that are studied under a specified policy.
It relies on the identification of a `behavior' that contributes most to the cost,
such as when the underlying stochastic dynamics (say, the suitably
normalized multidimensional queue length) lie close to a specific path.
In the case of controlled dynamics,
this path is formulated as the control selected by the maximizing player in the DG.
To obtain the lower bound one must consider an arbitrary sequence of policies, and
then the challenge stems from the fact that
different policies for the queueing model may give rise to different such paths.
A brute force approach of identifying an optimal path for each arbitrary policy
seems intractable.

The argument uses instead properties of the DG studied in \cite{ata-coh}.
It has been shown that this DG, specified in terms of multi-dimensional dynamics,
can be reduced to one dimension. The one-dimensional state corresponds to
the (suitably normalized) total workload in the system.
Moreover, there is a threshold, denoted by $\beta_0$, dictating
the behavior of both players. When the workload is below this threshold,
there is a certain fixed path that guarantees attaining at least the game's value
under any action of the minimizer (although it need not be optimal).
When the threshold is exceeded, there is no such fixed path.
However, the following fact can be used. As long as the workload remains above
the level $\beta_0$, the minimizer encounters an accumulated
loss, which is higher than the cost of an immediate rejection
to the level $\beta_0$. We identify a suitable path for the maximizer
that is effective until the time when the threshold is reached. We focus
on this path when workload is above $\beta_0$, and switch to the path alluded to
above, when it is below $\beta_0$.
However, this switching time depends on the policy,
and so it is random and varies with the scaling parameter.
Accordingly, the argument uses time discretization, where each
one of a finite collection of possible switching times is estimated separately.

The upper bound is obtained by constructing a policy for the QCP
for which the cost converges to the DG's value.
There is a naive way of interpreting the DG solution as a control policy for the QCP.
However, the two components of this policy, corresponding to rejection
and service effort allocation, impose contradictory requirements.
For rejections must occur only from a specific
class, and only when the total workload in the system exceeds the aforementioned
threshold. On the other hand, it stems from the game solution that
the service allocation policy must cause the (suitably normalized) multidimensional
queue length processes to evolve along a certain curve in state space,
denoted in this paper by $\gamma$. One can express the fact that buffers are finite
by requiring that these queueing
processes always lie in a certain hyper-rectangular domain, denoted by $\calX$;
the curve $\gamma$ happens to intersect the boundary of the domain $\calX$.
When the multidimensional queueing process is on (certain parts of) that boundary,
one of the buffers must be full, and then even small stochastic fluctuations require
rejections so as to meet the buffer size constraint. As a consequence,
rejections will not always adhere to the rejection policy alluded to above.

We address this issue by aiming at a curve $\gamma^a$, that approximates $\gamma$
but does not intersect the boundary of the domain. The main body
of work is to develop estimates that show that the queueing processes
evolve along this approximate curve, up to negligible probabilities,
and that as a result both elements of the policy are respected
with sufficiently high probability.

It is well known that an optimal control over a finite time interval
is generically of a {\it nonstationary} feedback type
(see, for example, Sections III.8--9 in \cite{Fleming2006}, where the nonstationary
optimal feedback is characterized by (8.5) and the stationary optimal feedback by (9.9),
corresponding to a problem set on a finite time horizon
and, respectively, an infinite time horizon with a discounted cost).
Despite that, it was anticipated in \cite{ata-coh},
and established in this paper, that for the setting studied here,
a RS cost set over a finite time horizon $[0,T]$
gives rise to a {\it stationary} feedback provided that $T$ is
sufficiently large. Indeed, the policy we present for the QCP has this
feature, which makes it simple as compared to policies based on
time-varying characteristics. We consider this as one of the main aspects
of this paper's contribution.

The organization of the paper is as follows. Section \ref{sec2} presents the model,
the MD scaling and the main result, which states that the MD QCP's value converges
to that of the DG. Section \ref{sec3} collects a few results from \cite{ata-coh}
required for the proof. Section \ref{sec_lower}
gives a lower bound on the QCP's value asymptotics in terms of
the DG's value, and Section \ref{sec_upper} finds a nearly optimal policy derived from
the game's optimal strategy. Together, Sections \ref{sec_lower} and \ref{sec_upper}
provide the proof of the main result.
Some auxiliary results appear in the Appendix.

\skp

We use the following notation.
For a positive integer $k$ and $a,b\in\R^k$, $a\cdot b$ denotes the usual
scalar product, while $\|\cdot\|$ denotes Euclidean norm.
$\{e_1,\ldots,e_k\}$ is the standard basis of $\R^k$.
We denote $[0, \iy)$ by $\R_+$.
For $0<T<\iy$ and
a function $f:\R_+\to\R^k$, $\|f\|_T=\sup_{[0,T]}\|f\|$,
while $\osc_T(\del,f)=\sup\{\|f(u)-f(t)\|:0\le u\le t\le (u+\del)\w T\}$.
Denote by $\AC([0,T],\R^k)$, $\calC([0,T],\R^k)$ and $\calD([0,T], \R^k)$ the
spaces of absolutely continuous functions [resp.,
continuous functions, functions that are
right-continuous with finite left limits (RCLL)] mapping $[0,T]\to\R^k$.
Write $\AC_0([0,T],\R^k)$ and $C_0([0,T],\R^k)$ for the subsets of the corresponding
function spaces, of functions that start at zero.
Endow the space $\calD([0,T],\R^k)$ with the usual Skorohod topology.
For a collection $x_i\in\R$ indexed by $i\in\{1,\ldots,I\}$ ($I$ being
a positive integer), $x$ denotes the vector $(x_i)$.
A similar convention holds for $\R$-valued random variables $X_i$
and stochastic process $\{X_i(t), t\in\R_+\}$, $i\in\{1,\ldots,I\}$,
where $X$ and $\{X(t), t\in\R_+\}$ denote the $\R^I$-valued random variable and process.

\section{Model and results}\label{sec2}
\beginsec
\subsection{Model description}\label{sec_multi_stoch_model}

We consider a model with $I$ customer classes and a single server.
A buffer with finite room is dedicated to each customer class, and
upon arrival, customers are queued in the corresponding buffers,
or rejected by the system administrator.
Within each class, customers are served at the order of arrival, where
the server may only serve the customer at the head of each line.
Processor sharing is allowed, and so the server
is capable of serving up to $I$ customers of distinct classes simultaneously.
The model is defined on a probability space $(\Om,\calF,\PP)$. Expectation with
respect to $\PP$ is denoted by $\EE$.
The parameters and processes we introduce depend on
an index $n\in\N$, serving as the scaling parameter.
Arrivals occur according to independent
renewal processes, and service times are independent and identically distributed
across each class.
Let $\I=\{1,2,\ldots, I\}$.
Let $\la^n_i >0, n\in\N, i\in\I$ be given parameters, representing the {\it reciprocal mean
inter-arrival times} of class-$i$ customers.
Let $\{\IA_i(l) : l\in\N\}_{i\in\calI}$ be $I$ sequences of positive i.i.d.\ random variables with mean $\E[\IA_i(1)]=1$ and
variance $\sig^2_{i,IA}=\mbox{Var}(\IA_i(1))\in (0,\infty)$. With $\sum_1^0=0$,
the number of arrivals of class-$i$ customers up to time $t$, for the $n$-th system,
is given by
\begin{align}\label{20}
A^n_i(t):=\sup\Big\{ l\geq 0\ : \sum_{k=1}^l\frac{\IA_i(k)}{\la^n_i}\leq t\Big\}, \quad t\geq 0.
\end{align}
For a collection $\xi_i$, $i\in\calI$ of stochastic processes we will always write
$\xi$ for $(\xi_i)_{i\in\calI}$. Thus, in particular,
$A^n$ is the $I$-dimensional process $(A^n_i)_{i\in\calI}$.

Similarly we consider another set of parameters $\mu^n_i>0, n\in\N, i\in\calI,$
representing {\it reciprocal mean service times}.
We are also given $I$ independent sequences $\{\ST_i(l) : l\in\N\}_{i\in\calI}$ of positive i.i.d.\ random variables (independent also of the sequences $\{\IA_i\}$)
with mean $\E[\ST_i(1)]=1$ and variance $\sig^2_{i,ST}=\mbox{Var}(\ST_i(1))\in (0,\infty)$.
The time required to complete the service of the $l$-th class-$i$ customer is given
by $\ST_i(l)/\mu^n_i$, and the {\it potential service time} processes are defined as
\begin{align}\notag
S^n_i(t):=\sup\Big\{ l\geq 0\ : \sum_{k=1}^l\frac{\ST_i(k)}{\mu^n_i}\leq t\Big\}, \quad t\geq 0.
\end{align}
We consider the {\it moderate deviations rate parameters} $\{b_n\}$, that form a sequence,
fixed throughout, with the property that $\lim b_n=\infty$ while $\lim b_n/\sqrt{n}=0$,
as $n\to\iy$.
The arrival and service parameters are assumed to satisfy the following conditions.
As $n\to\iy$,
\begin{equation}\label{11}
\frac{\la^n_i}{n}\to\la_i\in(0,\iy),
\quad
\frac{\mu^n_i}{n}\to\mu_i\in(0,\iy),
\end{equation}
\begin{equation}\label{12}
\tilde\la^n_i:=\frac{1}{b_n\sqrt{n}}(\la^n_i-n\la_i)\to\tilde\la_i\in(-\infty, \infty),
\quad
\tilde\mu^n_i:=\frac{1}{b_n\sqrt{n}}(\mu^n_i-n\mu_i)\to\tilde\mu_i\in(-\infty, \infty).
\end{equation}
The system is assumed to be critically loaded in the sense that the overall traffic intensity
equals 1, namely $\sum_1^I \rho_i=1$ where $\rho_i=\la_i/\mu_i$ for $i\in\calI$.

For $i\in\calI$, let $X^n_i$ be a process representing the number of class-$i$ customers in the $n$-th system. Denote the number of rejection of class-$i$ arrivals
until time $t$ in the $n$-th system by $R^n_i(t)$.
With $\mathbb{S}=\{x=(x_1,\ldots,x_I)\in\R_+^I:\sum x_i\le1\}$, let
$B^n$ be an $\mathbb{S}$-values process, whose $i$-th component represents
the fraction of effort devoted by the server to the class-$i$ customer at the head of the line.
Then the number of service completions of class-$i$ jobs during
the time interval $[0,t]$ is given by $S^n_i(T^n_i(t))$,
where
\begin{equation}\label{eq2007}
T^n_i(t):=\int_0^t B^n_i(u)du
\end{equation}
is the time devoted to class-$i$ customers by time $t$.
With an abuse of notation, we often write
$S^n\circ T^n$ for $(S^n_1\circ T^n_1,\ldots,S^n_I\circ T^n_I)$.
We have the balance equation
\begin{equation}\label{eq2008}
X^n_i(t)=X^n_i(0)+A^n_i(t)-S^n_i(T^n_i(t))-R^n_i(t).
\end{equation}
For simplicity, the initial conditions $X^n(0)=(X^n_1(0),\ldots,X^n_I(0))$
are assumed to be deterministic. We also assume
\begin{equation}\notag%\label{eq2017}
\tilde{X}^n(0)\to \bar x=(x_1,\ldots,x_I) \ \mbox{as}\ n\to\infty,\quad i\in\I.
\end{equation}
Note that, by construction, the arrival and potential service processes have RCLL paths,
and accordingly, so does $X^n$.

The MD-scaled version of the queue length process satisfies
\begin{equation}\label{10}
\tilde{X}^n(t):=\frac{1}{b_n\sqrt{n}}X^n(t)\in\calX:=\prod_{i=1}^I[0,D_i],
\qquad t\ge0,
\end{equation}
where $D_i>0$ are fixed constants. Thus the size of buffer $i$
is given by $b_n\sqrt n D_i$.
Additional MD-scaled processes are
\begin{align}\notag%\label{eq2012}
 \tilde{A}^n_i(t)&=\frac{1}{b_n\sqrt{n}}(A^n_i(t)-\la^n_i t),
  \quad \tilde{S}^n_i(t)=\frac{1}{b_n\sqrt{n}}(S^n_i(t)-\mu^n_i t),
  \quad \tilde{R}^n_i(t)=\frac{1}{b_n\sqrt{n}}R^n_i(t).
\end{align}
The process $U^n:=(B^n, R^n)$ is regarded as a control, that is determined based on
observations from the past events in the system.
The precise definition of an admissible control is as follows.
Given $n$, the processes $U^n$ and $X^n$ are said to be an {\it admissible control}
and the {\it corresponding queue length process}
if the sample paths of $U^n$ lie in $\calD([0,\iy),\mathbb{S})$, \eqref{10} holds,
and
\begin{itemize}
\item $U^n$ is adapted to the filtration
 $
 \s\{A^n_i(u),S^n_i(T^n_i(u)), i\in\I, u\leq t\},%\label{eq2009}
 $
 where $T^n$ is given by \eqref{eq2007};
\item For $i\in\I$ and $t\ge0$, one has
\begin{equation}
X^n_i(t)=0 \quad \text{implies} \quad B^n_i(t)=0. \label{eq2010}
\end{equation}
\end{itemize}
Denote the class of all admissible controls $U^n$ by $\Uf^n$.
Note that this class depends on $A^n$ and $S^n$, but we consider these processes
to be fixed. For $U^n\in\Uf^n$ and the corresponding queue length process $X^n$,
the processes $\tilde R^n$ and $\tilde X^n$ are referred as the scaled rejection and
queue length process corresponding to $U^n$.

Throughout, we assume the finite exponential moment condition, namely
\begin{assumption}
  \label{assn1}
  There exists $u_0>0$ such that for $i\in\I$,
  $\E[e^{u_0\IA_i}]$ and $\E[e^{u_0\ST_i}]$ are finite.
\end{assumption}
As shown in \cite{Puhal-Whitt}, under this condition,
the scaled processes $(\tilde A^n,\tilde S^n)$ satisfy a {\it moderate deviation
principle}. Namely, for $k=1,2$, let $\Jr_k(T,\cdot)$ be functions mapping $\calD([0,T],\R^I)$ to
$[0,\iy]$ given by
\begin{align}\label{eq2018}
\Jr_k(T,\psi)=\left\{\begin{array}{ll}
               \sum_{i=1}^I \cc_{i,k}
               \int_0^T\dot\psi_i(u)^2du &\ \mbox{if all}\ \psi_i\in\AC_0([0,T],\R),
\\
                \infty & \ \mbox{otherwise},
              \end{array}
\right.
k=1,2,
\end{align}
for $\psi=(\psi_1,\ldots,\psi_I)\in \calD([0,T],\R^I)$, where
\begin{align}\notag%\align{eq2019}
\cc_{i,1}= \frac{1}{2\la_i\sig^2_{i,IA}}\quad\text{and}\quad \cc_{i,2}=\frac{1}{2\mu_i\sig^2_{i,ST}},\quad i\in\calI.
\end{align}
Let $\Jr(T,\psi)=\Jr_1(T,\psi^1)+\Jr_2(T,\psi^2)$
for $\psi=(\psi^1,\psi^2)\in \calD([0,T],\R^{2I})$. Note that $\Jr$ is lower semicontinuous with compact level sets. Then one has
\begin{proposition}
[\cite{Puhal-Whitt}]
\label{moderate}
Under Assumption \ref{assn1}, for $T>0$ fixed, the following holds.
For every closed set $F\subset \calD([0,T],\R^{2I})$,
$$
\limsup\frac{1}{b_n^2}\log\p((\tilde{A}^n,\tilde{S}^n)\in F)
\leq -\inf_{\psi\in F}\Jr(T,\psi),
$$
and for every open set $G\subset \calD([0,T],\R^{2I})$,
$$
\liminf\frac{1}{b_n^2}\log\p((\tilde{A}^n,\tilde{S}^n)\in G)
\geq -\inf_{\psi\in G}\Jr(T,\psi).
$$
\end{proposition}

To present the RS control problem, fix $\uh,\ur\in(0,\infty)^I$.
Given $T\in(0,\iy)$ and $n$, the cost associated with a control $U^n\in\Uf^n$ is given by
\begin{equation}\notag%\label{eq2025}
J^n(T,\tilde X^n(0),U^n)=\frac{1}{b_n^2}
\log\E\Big[\int_0^T
e^{b^2_n[\int_0^t \uh\cdot\tilde{X}^n(u)du+\ur\cdot\tilde{R}^n(t))]}dt\Big],
\end{equation}
where $\tilde X^n$ and $\tilde R^n$ are the rescaled queue length and rejection processes
corresponding to $U^n$.
In the above notation, we have emphasized the dependence on the initial state $\tilde X^n(0)$.
For background and a motivating discussion about this type of cost
the reader is referred to \cite{Fleming2006b}. The value of interest is given by
\begin{equation}\notag%\label{eq2026}
V^n(T,\tilde X^n(0))=\inf_{U^n\in\Uf^n} J^n(T,\tilde X^n(0),U^n).%\label{eq020}
\end{equation}

\subsection{The differential game and main result}

Whereas the scaled queue length process $\tilde X^n$ is multidimensional, it is suggested
in \cite{ata-coh} that it is governed by a DG defined in terms of one-dimensional
dynamics. The main result of this paper is the proof of this claim.
Before presenting the formulation of this game, it is useful to draw attention to
several stochastic processes that are themselves one-dimensional because the structure
of the DG is closely related to them.
Let
\begin{equation}\notag %\label{eq2026dda}
\theta^n=\Big(\frac{n}{\mu^n_1},\ldots,\frac{n}{\mu^n_I}\Big),
\quad
\theta=\Big(\frac{1}{\mu_1},\ldots,\frac{1}{\mu_I}\Big),
\end{equation}
and note that $\theta^n\to\theta$, by \eqref{11}.
Denote $\theta_{\rm min}=\max_i\theta_i$ and $\theta_{\rm max}=\max_i\theta_i$.
Let also
\[
D^n=\sum_i\theta^n_iD_i,\quad D=\sum_i\theta_iD_i.
\]

It follows from the balance equation \eqref{eq2008} that
\begin{equation} \label{eq2013}
 \tilde{X}^n_i=\tilde{Y}^n_i+\tilde{A}^n_i-\tilde{S}^n_i\circ T^n_i
 + \tilde Z^n_i-\tilde{R}^n_i,
\end{equation}
where we denote
\begin{equation}\label{eq2014}
\tilde Z^n_i(t)=\frac{\mu^n_i}{n}\frac{\sqrt{n}}{b_n}(\rho_i t-T^n_i(t)),
\qquad \tilde{y}^n_i=\tilde\la^n_i-\rho_i\tilde\mu^n_i,
\qquad \tilde Y^n_i(t)=\tilde X^n_i(0)+\tilde y^n_it.
\end{equation}
Define
\begin{align}\label{16}\notag
 \check{X}^n &=\theta^n\cdot\tilde{X}^n ,
  \quad  \check{A}^n =\theta^n\cdot\tilde{A}^n ,
  \quad  \check{S}^n =\theta^n\cdot\tilde{S}^n ,\\
  \quad  \check{R}^n &=\theta^n\cdot\tilde{R}^n ,
  \quad  \check{Y}^n =\theta^n\cdot\tilde{Y}^n,
  \quad \check{Z}^n =\theta^n\cdot\tilde{Z}^n ,
  \quad \check{y}^n =\theta^n\cdot\tilde{y}^n .
\end{align}
Also, let
\begin{align*}
\tilde\bS^n(t_1,\ldots,t_I)&=(\tilde S^n_1(t_1),\ldots,\tilde S^n_I(t_I)),
\\
\check\bS^n(t_1,\ldots,t_I)&=\theta^n\cdot\tilde\bS^n(t_1,\ldots,t_I)
=\sum_i\theta^n_i\tilde S^n_i(t_i).
\end{align*}
Note that the sample paths of $\tilde S^n$ [resp., $\check S^n$, $\tilde\bS^n$, $\check\bS^n$]
map $\R_+\to\R^I$ [resp., $\R_+\to\R_+$, $\R_+^I\to\R^I$, $\R_+^I\to\R$].
Next, the process
\begin{equation}\label{eq2015}
\check Z^n=\sum_i\frac{n}{\mu^n_i}\tilde Z^n_i \quad \text{is nonnegative and nondecreasing,}
\end{equation}
thanks to the fact that $\sum_iB^n_i\le1$ while $\sum_i\rho_i=1$.
Also, for every $i\in\I$, the process $\tilde{R}^n_i$ is nondecreasing.
Thus by \eqref{eq2013},
\begin{equation}
  \label{02}
  \check X^n=\check Y^n+\check A^n-\check \bS^n\circ T^n+\check Z^n-\check R^n,
\end{equation}
where $\check Z^n$ and $\check R^n$ are nonnegative, nondecreasing processes.
Moreover,
\begin{equation}\label{eq2026c1}
 \check{X}^n(t)\in[0,D^n], \qquad t\ge0.
\end{equation}

It follows from the contraction principle that $(\check A^n(t),\check S^n(t))$, $t\in[0,T]$,
satisfy the MDP with the rate function
$\Ir(T,\psi)=\Ir_1(T,\psi^1)+\Ir_2(T,\psi^2)$,
$\psi=(\psi^1,\psi^2)\in\calD([0,T],\R^2)$,
where
\begin{equation}\label{eq2026f}
\Ir_k(T,\psi^k)=\left\{\begin{array}{ll}
               \cc_k
               \int_0^T(\dot\psi^k)^2(u)du &\ \mbox{if}\ \psi^k\in\AC_0([0,T],\R),
\\
                \infty & \ \mbox{otherwise},
              \end{array}
\right.
\end{equation}
\begin{equation}\label{eq2026ff}
\cc_1:=\left(\sum_{i=1}^{I}\frac{2\rho_i\sig^2_{i,IA}}{\mu_i}\right)^{-1},\quad\text{and}\quad \cc_2:=\left(\sum_{i=1}^{I}\frac{2\rho_i\sig^2_{i,ST}}{\mu_i}\right)^{-1}.
\end{equation}
See Lemma \ref{lem_underline_psi}.

To define the DG,
Denote $x=\theta\cdot\bar x$,
$y:=\lim\check y^n=\sum_i\theta_i(\tilde\la_i-\rho_i\tilde\mu_i)$, and
\[
\by(t)=x+yt,\qquad t\in\R_+,
\]
and let
\begin{equation}\label{eq2027}
\calP=\calC_0([0,\iy),\R),
\quad
\calE=\{\xi\in \calD([0,\iy),\R_+):\xi \text{ is nondecreasing}\}.
\end{equation}
Endow both spaces with the topology of uniform convergence on compacts.
Given $\psi=(\psi^1,\psi^2)\in \calP^2$ and $(\zeta,\vr) \in \calE^2$, the {\it dynamics associated
with the initial condition $x$ and the data $\psi,\zeta,\vr$} is defined as
\begin{equation}\label{eq2028}
\ph= \by+ \psi^1-\psi^2+\zeta-\vr.
\end{equation}
The game is played by a maximizing player that selects $\psi=(\psi^1,\psi^2)$
and a minimizing player that selects $(\zeta,\vr)$.
We sometimes write the dependence of the dynamics on the data as
$\ph[x,\psi,(\zeta,\vr)]$.
There is an analogy between the above equation and equation \eqref{02},
and between the condition that $\zeta$ and $\vr$ are
nondecreasing and property \eqref{eq2015}. The control $\zeta$
stands for the scaled idle time process $\tilde{Z}^n$ and $\vr$
stands for the scaled rejection process $\tilde{R}^n$.
The following condition, analogous to property \eqref{eq2026c1}, will also be required,
namely
\begin{equation}\label{eq2029}
\ph(t)\in[0,D],\qquad t\ge 0,
\end{equation}
where
\[
D=\sum_{i=1}^I\theta_iD_i=\limn D^n.
\]
A measurable mapping $\alpha:\calP^2\to \calE^2$ is called a {\it strategy for
the minimizing player} if it satisfies the causality property:
for every $\psi, \tilde \psi \in \calP^2$
and $t\in[0,\iy)$,
\begin{equation}\label{eq2030}
\text{$\psi(u)=\tilde\psi(u)$
for every $u\in[0,t]$\quad implies\quad $\alpha[\psi](u)=\alpha[\tilde\psi](u)$
for every $u\in[0,t]$.}
\end{equation}
Given an initial condition $x$,
a strategy $\al$ is said to be {\it admissible for the initial condition $x$}
if, whenever $(\psi^1,\psi^2)\in \calP^2$ and
$(\zeta,\vr)=\al[\psi]$, the corresponding dynamics \eqref{eq2028} satisfies the buffer
constraint \eqref{eq2029}. We denote by $A_x$ the collection of admissible
strategies for the initial condition $x$.

We now describe the components of the cost function.
For $w\in\R_+$, denote
\begin{equation}\label{eq2026g}
 h(w)=\inf\{\uh\cdot \xi : \xi\in\calX,\ \theta\cdot \xi=w\}.
\end{equation}
By the convexity of the set $\calX$, $h$ is convex. Moreover, $ h(w)\ge0$ for $w\ge 0$ and equality holds if and only if $w=0$. Therefore, $ h$ is
strictly increasing on $[0,D]$. Let
\begin{equation}\notag%\label{eq20237}
 r=\min\{\ur\cdot \xi: \xi\in\R_+^I, \utheta\cdot \xi=1\}.
\end{equation}
It is easy to see that
\begin{equation}\label{eq20238aa}
 r=\min\{r_i\mu_i:i\in\I\}=r_{i^*}\mu_{i^*},
\end{equation}
where $i^*$ is an index (fixed throughout) that minimizes $r_i\mu_i$.

The index $i^*$ indicates the class that has lowest rejection cost per unit of workload.
It plays an important role in Section \ref{sec_upper} where our AO policy is presented;
specifically, the policy is aimed at rejecting jobs from this class only.

Given $x\in[0,D]$, $T\in\R_+$, $\psi=(\psi^1,\psi^2)\in \calP^2$, and $(\zeta,\vr)\in \calE^2$, we define the cost until time $T$ by
\begin{equation}\label{01}
c(x,T,\psi,\zeta,\vr)=\int_0^Th(\ph(t))dt+r\vr(T)-\Ir(T,\psi),%\label{025}
\end{equation}
where $\ph$ is the corresponding dynamics.
The value of the game is defined by
\begin{equation}\label{04}
V(x)=\inf_{\alpha\in A_x}\sup_{\psi\in \calP^2, T\in\R_+}c(x,T,\psi, \alpha[\psi]).%\label{026}
\end{equation}
We call $\psi$ the \emph{path control} and the $T$ a \emph{time control}, or sometimes
the {\it termination time}. Note that both are controlled by the maximizing player.

We sometimes use the notation $V_{h,r}$ for $V$ when we want to emphasize
the dependence on the function $h$ and the constant $r$. Namely, $V_{\hat h,\hat r}(x)$
is defined as $V(x)$ with $(\hat h,\hat r)$ in place of $(h,r)$ in \eqref{01}.

Recall that $\limn \tilde X^n(0) =x$. Our first main result is the following.
\begin{theorem}\label{main}
Let Assumption \ref{assn1} hold. Then
for all sufficiently large $T$,
\[\lim_{n\to\iy} V^n(T,\tilde X^n(0))=V(x).
\]
\end{theorem}
The second main result of this paper is Theorem \ref{thm_upper}, that constructs a policy
for the QCP, which is AO.

\section{Some useful properties of the game}\label{sec3}\beginsec

We briefly mention some results from \cite{ata-coh} regarding the DG,
to be used in the sequel.
Set $\cc=(\cc_1^{-1}+\cc_2^{-1})^{-1}$. The following proposition follows by Lemma 3.1, Proposition 3.1, and Theorems 3.1, 3.2, and 3.5 in  \cite{ata-coh}. The contribution of the latter is to deduce part (i) below.
\begin{proposition}\label{prop_game_value}
(i) For $T\in\R_+$, set
\[
V(T,x)=\inf_{\alpha\in A_x}\sup_{\psi\in \calP^2, t\in[0,T]}c(x,t,\psi, \alpha[\psi])
\]
(compare with \eqref{04}). Then $V(T,x)=V(x)$ for all sufficiently large $T$.
\\
(ii) If $-y<r/(4\cc)$ then for every $x\in[0,D]$ one has $V(x)=\iy$.\\
\noi (iii) If $-y\ge r/(4\cc)$ then
\begin{equation}\label{eq2041}
V(x)=\left\{\begin{array}{ll}
               \displaystyle\int_0^x 2\cc\Big(-y-\sqrt{y^2-\frac{h(u)}{\cc}}\Big)du,
               &\ 0\leq x\leq \Barr,
\\ \\
                V(\Barr)+r(x-\Barr), & \  \Barr<x\leq D,
              \end{array}
\right.
\end{equation}
where,\footnote{In case $\frac{r^2}{4\cc}+ry\ge-h(0)=0$ we get by (i) above that $V(x)=\iy$ and we do not define $\Barr$.}
\begin{equation}\label{eq2042}
  \beta_0=\begin{cases}
    \ds h^{-1}\Big(\frac{-r^2}{4\cc}-ry\Big), & \ds -h(D)\le \frac{r^2}{4\cc}+ry\le-h(0),\\ \\
    D, & \ds\frac{r^2}{4\cc}+ry<-h(D).
  \end{cases}
\end{equation}
\end{proposition}

We now present an optimal strategy for the minimizer. This strategy plays an important role in proving the upper bound and in finding an AO policy in the multidimensional stochastic problem. The minimizer's optimal strategy is of a $\beta$-barrier form.
Informally, this is a strategy that uses the minimal control $(\zeta,\vr)$
so as to keep the dynamics $\ph$ in $[0,\beta]$ at all times.
In the definition that follows, and throughout the paper, we denote the
Skorohod map on an interval $[a,b]$ by $\Gam_{[a,b]}$; see Appendix \ref{sec_skorohod}.

\begin{definition}\label{def_barrier}
Fix $(x,\beta)\in[0,D]^2$. The strategy $\al_\beta=(\al_{\beta,1},\al_{\beta,2})$
is called a $\beta$-barrier strategy if for every $\psi\in\calP^2$ one has
$(\ph,\al_{\beta,1},\al_{\beta,2})[\psi]=\Gam_{[0,\beta]}(\psi)$.
\end{definition}

The next proposition follows by Proposition 3.1 and Theorems 3.1, 3.2, and 3.4 in \cite{ata-coh}.
\begin{proposition}\label{prop_game_optimal_strategy}
The $\Barr$-barrier strategy, $\al_\Barr$, is an optimal strategy.
\end{proposition}

We provide two propositions that are useful in the proof of the lower bound. For this we present two path controls, $\psi^*$ and $\psi^\sharp_x$ associated with an initial state in the intervals $(\Barr,D]$ and $[0,\Barr)$ respectively.

Fix $x\in(\Barr,D]$. Let $\trian>0$ be such that $x>\Barr+\trian$.
Fix $(\zeta,\vr)\in\calP^2$.
Define
\begin{equation}\notag%\label{eq2043}
\psi^\sharp(t)=(rt/(2\cco),-rt/(2\cct)),\qquad 0\le t\le\tau_\Del,
\end{equation}
where
\begin{equation}\notag%\label{eq2044}
\tau_\trian=\tau_{(\zeta,\vr),\trian}:= \inf\{t\ge 0 : \ph[x,\psi^\sharp,(\zeta,\vr)]\leq\Barr+\trian \}
\end{equation}
is the first time that the dynamics, $\ph:=\ph[x,\psi^\sharp,(\zeta,\vr)]$, cross $\Barr+\trian$.
The following proposition is Proposition 3.2 in \cite{ata-coh}.
\begin{proposition}\label{prop_above_barrier}
For every $(\zeta,\vr)\in\calP^2$ such that $\vr(0)-\zeta(0)<x-(\Barr+\trian)$ one has
\begin{align}\label{eq2045}
\int_0^{\tau_{\trian}}h(\ph(t))dt+r\vr(\tau_{\trian})-\Ir(\tau_{\trian},\psi^\sharp)> r(x-(\Barr+\trian)).
\end{align}
\end{proposition}
Note that the l.h.s.~is the cost associated with $x,\psi^\sharp$, and $(\zeta,\vr)$ incurred until the time the dynamics cross $\Barr+\trian$, whereas the r.h.s.~gives the cost of an immediate rejection of $x-(\Barr+\trian)$.
The result thus implies that if $x>\Barr$ then the minimizer will reject $x-\Barr$ units of mass
at time zero.

Next, fix $x\in [0,\Barr)$. Let
\begin{align}\label{13}
\psi^*_x=\left(\frac{\cc}{\cco}\om^*_x,\frac{-\cc}{\cct}\om^*_x\right),
\end{align}
where $\om^*_x\in \calC([0,\tau_x^*),\R)$ is the unique solution of
\begin{equation}\label{14}
\dot\om^*_x(t)=\frac{\dot{V}(x+yt+ \om^*_x(t))}{2\cc},\qquad  t\ge 0,
\end{equation}
with $\om^*_x(0)=0$ and
\begin{align}\label{102}
\tau_x^*= \int_0^{x} 1/\sqrt{y^2-h(\xi)/\cc}\;d\xi.
\end{align}
Existence and uniqueness for $\om^*_x$ over the time interval $[0,\tau_x^*]$ are shown in Section 3.6.2 in \cite{ata-coh}.

The next proposition, which follows by Proposition 3.3, Theorem 3.5, and equation (87) in \cite{ata-coh}, states that by using the path control $\psi^*_x$ and by choosing the time control to be the first time that the actual dynamics of the game hit zero the maximizer can guarantee that the cost will be at least $V(x)$. The proposition is valid since the function $h$ is convex.
\begin{proposition}\label{prop_below_barrier}
Fix $x\in[0,\Barr)$.
For every $\vr\in\calP$ one has
\begin{equation}\label{eq2047b}
c(x,\tau_x,\psi^*_x, (0,\vr))\ge V(x),
\end{equation}
where $\tau_x:=\tau[x,\psi^*_x,(0,\vr)]$ is the first time that the dynamics $\ph[x,\psi^*_x,(0,\vr)]$ hits zero.
Moreover,
\begin{equation}\label{eq2047c}
\tau_x\le \tau_x^*= \inf\{t\ge 0 :\by(t)+\psi^{*,1}_x(t)-\psi^{*,2}_x(t)=0 \}.
\end{equation}
\end{proposition}

\section{Lower bound}\label{sec_lower}
\beginsec

Recall that from Proposition \ref{prop_game_value}.(i), for sufficiently large $T$, $V(T,x)=V(x)$. Hence, now onwards we compare the value function of the QCP to $V(x)$.
\begin{theorem}\label{thm_lower}
Let Assumption \ref{assn1} hold. Then for $T$ sufficiently large,
\[
\liminf_{n\to\iy} V^n(T,\tilde X^n(0))\geq  V(x).
\]
\end{theorem}

We present two lemmas that together yield Theorem \ref{thm_lower}.
The first provides a lower bound on the RS cost
for an arbitrary sequence of policies in terms of an expression involving only
the one-dimensional processes. The latter is further bounded by the DG value function, in the
second lemma.
\begin{lemma}\label{prop_lower1}
Fix a sequence of admissible controls $\{U^n\in\Uf^n\}$, $n\in\N$. Then
for every $T>0$, $\delta\in(0,T)$ and $\eps>0$ one has
\begin{align}\notag%\label{eq2102}
\liminf_{n\to\iy} J^n (T,\tilde X^n(0),U^n) \ge
\liminf_{n\to\iy}\frac{1}{b_n^2}\log\EE\left[e^{b_n^2\left(\int_0^{T-\del} h(\check{X}^n(u))du+r(1-\eps)\check{R}^n(T-\del)\right)} \right]-\eps.
\end{align}
\end{lemma}

\begin{lemma}\label{prop_lower2}
Fix $\{U^n\}$ as in Lemma \ref{prop_lower1}.
Then there exists $\T>0$ such that for every $0<\eps<1/2$ one has
\begin{align}\label{eq2103}
\liminf_{n\to\iy}\frac{1}{b_n^2}\log\EE\left[e^{b_n^2\left(\int_0^{\T} h(\check{X}^n(u))du+r(1-\eps)\check{R}^n(\T)\right)} \right]\ge V_{h,r(1-\eps)}(x).
\end{align}
\end{lemma}

\noi{\bf Proof of Theorem \ref{thm_lower}:}
Combining Lemma \ref{prop_lower1} and Lemma \ref{prop_lower2},
for any $T$ and $\del>0$ such that $T-\del>\bar T$, and any $\eps\in(0,1/2)$,
\begin{align}\notag%\label{eq2104}
\liminf_{n\to\iy}J^n(T,\tilde X^n(0),U^n)\ge\liminf_{n\to\iy} J^n (\bar T,\tilde X^n(0),U^n) \ge
V_{h,r(1-\eps)}(x)-\eps.
\end{align}
 From \eqref{eq2041}--\eqref{eq2042} it follows that
\begin{align}\notag%\label{eq2105}
\lim_{\eps\to 0} V_{h,r(1-\eps)}(x) =V_{h,r}(x)=V(x),
\end{align}
and the result follows.
\qed

\skp\noi{\bf Proof of Lemma \ref{prop_lower1}:}
Fix $\{U^n\}$, $T>0$, $\del>0$ and $\eps>0$. Then
\begin{equation}\label{eq2106}
\int_0^Te^{b_n^2\left(\int_0^t \uh\cdot\tilde X^n(u)du+\ur\cdot\tilde R^n(t)\right)}dt
\ge e^{b_n^2\left(\int_0^{T-\del} \uh\cdot\tilde X^n(u)du+\ur\cdot\tilde R^n(T-\del)\right)}
\del,
\end{equation}
where we used monotonicity of the integrand with respect to $t$. Next,
by the definition of $h$ and $r$,
\begin{align}\label{eq2107a}
\uh\cdot\tilde X^n\ge h(\utheta\cdot\tilde X^n)
\quad\text{and}\quad
\ur\cdot\tilde R^n\ge r\utheta\cdot\tilde R^n.
\end{align}
Since $\utheta^n\to\utheta$, $\tilde X^n$ takes values in a fixed, compact set,
and $h$ is uniformly continuous on this set, it follows that
for sufficiently large $n$,
\begin{align}\label{eq2107d}
\int_0^T h(\utheta\cdot\tilde X^n(u))du\ge \int_0^T h(\check{X}^n(u))du-\eps\quad\text{and}\quad \theta_i\ge \theta^n_i(1-\eps),\;i\in\calI.
\end{align}
Combining \eqref{eq2106}, \eqref{eq2107a} and \eqref{eq2107d} yields the result.
\hfill $\Box$

\skp

In the rest of this section we prove Lemma \ref{prop_lower2}.

\skp\noi{\bf Proof of Lemma \ref{prop_lower2}:} Rather than working with general $\eps
\in(0,1/2)$, we consider a general $r>0$. For every $r>0$ we find $\T=\T(r)$
(in \eqref{eq2115}) that satisfies \eqref{eq2103} with $\eps=0$. As we will see, $\T$ is
continuous w.r.t.~$r$ and finite for $r>0$. Hence, we obtain that \eqref{eq2103}
is valid with a fixed $\bar T$ and all $\eps\in(0,1/2)$.
We thus turn to proving
\begin{align}\label{eq2107f}
\liminf\frac{1}{b_n^2}\log\EE\left[e^{b_n^2\left(\int_0^{\T} h(\check{X}^n(u))du+r\check{R}^n(\T)\right)} \right]\ge V(x).
\end{align}

\skp\noi{\bf Sketch of the proof:}
The lemma relates the stochastic control problem to the DG.
The strategies in the game are analogues of the policies in the stochastic problem,
whereas the controls selected by the maximizing player play a similar role
to the variational problem in Varadhan's lemma (Theorem 4.3.1 of \cite{DZ}).
Following the spirit of the proof of Varadhan's lemma, one focuses on the event that
the paths $(\tilde A^n,\tilde S^n)$, projected in the $\theta^n$ direction,
are in a neighborhood of a specific $\calP^2$-path control. The latter
is referred to as the reference path. One then shows that the process $\check{X}^n$
is in a neighborhood of the game dynamics obtained when the reference path is selected
by the maximizing player.

We now describe the reference path.
First, recall the balance equation \eqref{02}.
Consider paths $(\tilde A^n,\tilde S^n)$ such that $(\check{A}^n,\check{\bS}^n\circ T^n)$
are close to some path $\psi=(\psi^1,\psi^2)\in\calP^2$. Then $\check{X}^n$
is close to
$$
\by +\psi^{1}-\psi^{2} + \check{Z}^n -\check{R}^n.
$$
By the nonnegativity of $\check{Z}^n$, we have that $\check{X}^n$
is bounded from below by
$$
\ph^n= \by +\psi^{1}-\psi^{2}-\check{R}^n,
$$
up to a small error term.
The proof proceeds by comparing the process $\ph^n$
to the game dynamics with initial state $x$,
maximizer's path control $\psi$, and minimizer's control given by $(0,\check R^n)$.
This process cannot be regarded game dynamics because of the stochasticity
of the minimizer's control term, and the fact that it is not
attained by a strategy in the sense of the game. Moreover,
it is not guaranteed that $\ph^n$ takes values in $[0,D]$.

However, these obstacles can be treated.
Let us first describe the case $x\le \beta_0$, where the treatment is least
complicated.
In this case we consider the event, denoted by $O^n$, that
the paths $(\tilde A^n,\tilde S^n)$ are close to the path
$\psi^*_x$ identified in Proposition \ref{prop_below_barrier}
(in this rough sketch we do not quantify the term ``close'').
We then obtain
\begin{align}\notag
&\liminf\frac{1}{b_n^2}\log\EE\left[e^{b_n^2\left(\int_0^{\T} h(\check{X}^n(u))du+r\check{R}^n(\T)\right)} \right]\\\notag
&\qquad\ge \liminf\frac{1}{b_n^2}\log\EE\left[e^{b_n^2\left(\int_0^{\T} h(\check{X}^n(u))du+r\check{R}^n(\T)\right)} 1_{O^n}\right]\\\notag
&\qquad\ge c(x,\tau_x,\psi^*_x, (0,\vr))\ge V(x),
\end{align}
where the second inequality follows using the closeness of the data to $\psi^*_x$
and the structure of the function $c$, while the last inequality follows
from Proposition \ref{prop_below_barrier}.
Thus, Proposition \ref{prop_below_barrier}
allows us to focus on data that is close to one fixed path, $\psi^*_x$, for every $n\in\N$.

The situation is more subtle in the case
where $x>\beta_0$, as we do not have an analogue of Proposition \ref{prop_below_barrier}.
That is, we are unable to identify a single path that guarantees $V(x)$ as a lower bound on the cost
under all strategies. Proposition \ref{prop_above_barrier} proposes
how the maximizer in the game should act until the threshold
$\beta_0$ is reached, namely to use the control $\psi^\sharp$ (note that the time
when the threshold is reached depends on the minimizier's strategy).
Therefore we focus on data $(\tilde A^n,\tilde S^n)$ close to $\psi^\sharp$
until the workload in the stochastic model hits $\beta_0$.
The proposition then guarantees that up to this time the cost incurred is
bounded below by $r(x-\beta_0)$.
Once $\beta_0$ is reached, $\psi_{\beta_0}$ is used as explained above,
and \eqref{eq2041} is used to obtain $V(x)$ as a lower bound.

As we already mentioned, in the game,
the time when one switches from $\psi^\sharp$ to $\psi^*_{\beta_0}$
depends on the strategy. As far as the QCP is concerned, this means that one has
no control over the switching time, which may be random and vary with $n$.
The argument therefore uses time discretization (see Lemma \ref{lem_n1}),
with which each one of a finite collection of switching times is estimated separately.

\skp\noi{\bf Proof in details:}
We first prove the lower bound for the case $x>\Barr$. The other case, which
is simpler, is addressed at the end of the proof.
Fix $x>\Barr$ and $\eps_1>0$. Let $\trian>0$ be small enough so that
\begin{align}\notag%\label{eq2113
x>\Barr+\trian\qquad\text{and}\qquad r\trian+V(\Barr)-V(\Barr-\trian)\le \eps_1.
\end{align}
Such a $\trian$ exists since $V$ is continuous, see \eqref{eq2041}.
Let
\begin{align}\notag%\label{eq2114}
T_1&=\frac{V(x)+3+r(2+D-x)}{h(\Barr+\trian/2)-h(\Barr)},\\\notag%\label{eq310}
T_2&=\int_0^{\Barr-\trian} 1/\sqrt{y^2-h(\xi)/\cc}\;d\xi +1,
\end{align}
\begin{align}\label{eq2115}
\T=T_1+T_2.
\end{align}
Above, $T_1$ is obtained by considerations along these lines.
Start with \eqref{eq2143}, and consider only those policies
for which $V(x)+1>b_n^{-2}\log\E\exp\{b_n^2(\int_0^{T_1}h(\check X^n(u))du+r\check R^n(T_1))\}$.
Obtain a further lower bound on $V(x)+1$
by putting on the RHS the indicator of a certain event that assures
$\int_0^{T_1}h(\check X^n(u))du\ge h(\beta_0+\Delta/2)T_1$. Moreover, we show that
$\check R^n(T_1)\ge x-D-2+(y+\frac{r}{2s})T_1$.
When combined together with the definition of $\beta_0$ in \eqref{eq2042}, we get the above formula for $T_1$.
The formula for $T_2$ is based on the termination time of the game \eqref{102}.

The proof proceeds in several steps. In Step 1 we analyze
the process $T^n$. In Step 2, we consider the process
\begin{align}\label{101}
\ph^n_1: = \by +\psi^{\sharp,1}-\psi^{\sharp,2}-\check{R}^n
\end{align}
and the time $\tau^n_1$ when this process first hits $\Barr+\trian$.
We show that for any policy that is nearly optimal $\tau^n_1<T_1$ with
probability that is significant in the MD scale.

Next, in Step 3, we modify the construction of $\ph^n_1$ on the time interval
$[\tau^n_1,\iy)$, and in Step 4 show that the process thus constructed
hits zero before time $\bar T$, with probability 1.
In Step 5, we combine these results to obtain a lower bound on the
cost. Finally, in Step 6 we we take limits and obtain the result.

\skp\noi
{\bf Step 1 (Limit property of $T^n$):}
W.l.o.g.~assume that the sequence of policies $\{U^n\}$ satisfies
\begin{align}\label{eq2119}
V(x)+1>\frac{1}{b_n^2}\log\EE\left[e^{b_n^2\left(\int_0^{\T} \uh\cdot \tilde X^n(u)du+\ur\cdot\tilde R^n(\T)\right)} \right].
\end{align}
Denote
\[
\brho(t)=\rho t,\qquad t\in\R_+.
\]
\begin{lemma}\label{lem_k_m}
For every $m>0$ there exists $K>0$ such that for every $i\in\calI$ one has
\begin{equation}\label{05}
\limsup\frac{1}{b_n^2}\log \PP\left(\|T^n_i-\brho_i\|_{\T}
\geq \frac{b_n}{\sqrt{n}}K\right)\leq -m.
\end{equation}
\end{lemma}
\noi {\bf Proof:}\ Fix $i\in\calI$.
By \eqref{eq2014}, the l.h.s.\ of \eqref{05} equals
$\limsup\frac{1}{b_n^2}\log \PP\left(\frac{n}{\mu_i^n}\|\tilde Z^n_i\|_{\T} \geq K\right)$.
Using \eqref{eq2013} and the fact that $\mu_n/n\to\mu$, it suffices to prove
that for every $m>0$ there exists $K$ such that
\begin{equation}\label{eq2121}
\limsup\frac{1}{b_n^2}\log \PP\left(\|L^n\|_{\T} \geq K\right)\leq -m,
\end{equation}
for $L^n=\tilde X^n_i,\tilde A^n_i,\tilde S^n_i\circ T^n_i$ and $\tilde R^n_i$.
As far as $\tilde X^n_i$ is concerned, the above is immediate because the process is bounded.
For $L^n =\tilde A^n_i$, this property follows
from Proposition \ref{moderate} and from the fact
that for sufficiently large $K$ one has
\begin{align}\label{eq2123}
\inf\left\{\Jr_1(\T,\psi) : \psi\in\calP^I\quad\text{and}\quad \|\psi\cdot\ue_i\|_{\T}
\geq K\right\}= \Jr_1(\T,(Kt/\T)\ue_i)=\frac{1}{2\mu\sig^2_{IA}}\cdot\frac{K^2}{\T}.
\end{align}
Since the time change $T^n_i(t)\le t$ for all $t$, a similar conclusion holds
for $\tilde S^n_i\circ T^n_i$. Finally, for $\tilde R^n_i$,
note that by \eqref{eq2119},
\begin{align}\notag%\label{eq2124}
V(x)+1
&\geq
\frac{1}{b_n^2}\log\E\Big[ e^{b^2_n r_i\tilde R^n_i({\T})}\Big].
\end{align}
Hence by the Chebyshev's inequality,
$\PP(\tilde R^n_i ({\T})\geq K)\leq e^{-b_n^2(K +V(x)+1)}$.
The result follows since, by the monotonicity, $\tilde R^n_i(\T) =
\|\tilde R^n_i\|_{\T}$.
\hfill $\Box$

\skp

\noi{\bf Step 2 (Estimate on the time $\tau^n_1$):} We introduce some notation
that will be needed in the remainder of the proof.
Since $h$ is uniformly continuous on $[0,D]$, one can find $\delta_1>0$ such that
\begin{align}\label{eq2125b}
\osc_D(2\delta_1,h)<\eps_1.
\end{align}
One may take $\delta_1$ so that
\begin{align}\label{eq2126}
\delta_1<\min\{\trian/4,1\}.
\end{align}
Fix $m> \Ir({\T},\psi^*)+\Ir({\T},\psi^\sharp)+1+6\eps_1$.
Define the event
\begin{align}\notag%\label{eq2128}
E^n=E^n(K)=\{\| T^n_i-\brho_i\|_{\T} < \frac{b_n}{\sqrt{n}}K ,\text{for all } i\in\calI\}.
\end{align}
Using Lemma \ref{lem_k_m}, fix $K>0$ such that for all sufficiently large $n$
\begin{align}\label{06}
\frac{1}{b_n^2}\log \PP\left((E^n)^c\right)\leq -m.%\label{eq156}
\end{align}
From Lemma \ref{lem_underline_psi} it follows that there exists $\upsi^\sharp=(\upsi^{\sharp,1},\upsi^{\sharp,2})\in\calP^{2I} $ such that
\begin{align}\label{eq2128b}
(\utheta\cdot\upsi^{\sharp,1},\utheta\cdot\upsi^{\sharp,2}\circ\brho) =
(\psi^{\sharp,1},\psi^{\sharp,2})
\end{align}
and
\begin{align}\label{eq2128c}
\Jr_k(T,\upsi^{\sharp,k})= \Ir_k(T,\psi^{\sharp,k}) ,\quad k=1,2.
\end{align}
Note that for all large $n$,
\begin{equation}\label{eq2129}
\osc_T(\frac{b_n}{\sqrt{n}}K,\psi^{\sharp,2}_i)<\delta_2,\quad i\in\calI.
\end{equation}

For $\psi=(\psi^1,\psi^2)\in \calP^{2I}$, and $0<\delta, t<\infty$,
let
\begin{equation}\label{eq2130}
\A_{\delta,t}(\psi)=\{\upsi\in \calD([0,\bar T], \R^{2I}) \ :\ \|\upsi-\psi\|_t<\delta\},
\end{equation}
and
\begin{equation}\label{eq2134}
\Om_{\delta,t}^n(\psi)=\{(\tilde A^n,\tilde S^n)\in\A_{\delta,t}(\psi)\}.
\end{equation}

Recall from \eqref{101} the definition of
$\ph^n_1$, and let
\begin{equation}\notag%\label{eq2136}
\tau^n_1=\inf\{t\ge0 : \ph^n_1(t)\le \Barr+\trian\}.
\end{equation}
Divide the time interval $[0,\T]$ into $(\T/\nu)\in\N$ intervals of size
$\nu$ where\footnote{We use the convention that $\trian/ 0=\iy$.}
\begin{align}\label{eq2118}
\nu\le \min\Big\{\frac{4\cc\eps_1}{r^2}, \frac{\eps_1}{\cc y^2}, \frac{\trian}{|r/(2\cc)+y|}\Big\}.
\end{align}
Denote the intervals by $N_j=N_j(\nu)=[\nu j,\nu(j+1))$.
For every $n$, we define
an index $0\le j^n_1\le \lfloor\T_1/\nu \rfloor$
in such a way that we can estimate, from below, the probability
that the time $\tau^n_1$ belongs to the interval $N_{j_1^n}$. Let
\begin{align}\notag%\label{eq2137}
j^n_1=\underset{j\in\{0,\ldots,\lfloor T_1/\nu\rfloor\}}{\arg\max} \PP\left( \tau^n_1\in N_j
\mid \Om^n_{\delta_3,{T_1}}(\upsi^\sharp)  \right),
\end{align}
where $\del_2=\del_1/(8\theta_{\rm max}\sqrt{I})$.

\begin{lemma}\label{lem_n1}
Fix $\psi\in\calD([0,\bar T],\R^I)$ such that $\psi=\bar\psi^\sharp$ on $[0,(j^n_1+1)\nu]$.
Then
\begin{align}\label{eq2138}
\liminf_{n\to\iy}\PP\left( \tau^n_1\in N_{j^n_1}
\mid \Om^n_{\delta_3,{(j^n_1+1)\nu}}(\psi)\right) \ge -2\eps_1.
\end{align}
\end{lemma}
The proof of this lemma is differed to the end of the section.

\skp\noi{\bf Step 3 (Constructing a path beyond time $\tau^n_1$):}
In Lemma \ref{lem_n1} we focused on data
for which the process $(\check{A}^n,\check{\bS}^n\circ T^n)$ is near $\psi^\sharp$.
Now we will consider data for which this process is near
$\psi^\sharp$ up to time $(j^n_1+1)\nu$, and from that time on,
near $\psi^*_{\Barr-\trian}$.
Thus we focus on the reference path,
\begin{equation}\notag%\label{eq2149}
\psi^{\hs}(t)=\psi^{n,\hs}(t):=\left\{\begin{array}{ll}
               \psi^\sharp(t)  &\ 0\le t \le (j^n_1+1)\nu,
\\
                \psi^\sharp((j^n_1+1)\nu)+\psi^*_{\Barr-\trian}(t - (j^n_1+1)\nu) & \ t>(j^n_1+1)\nu.
              \end{array}
\right.
\end{equation}
Recall that $\tau^n_1$ is defined as the first time when $\ph^n_1\le\Barr+\trian$. Since $\by+\psi^{\sharp,1}-\psi^{\sharp,2}$ is continuous, and the jumps of $\check{R}^n$
are of size $(b_n\sqrt{n})^{-1}$, one has for sufficiently large $n$,
$\Barr< \ph^n_1(\tau^n_1)\le \Barr+\trian$. Moreover, since on the interval
$(\tau^n_1,(j^n_1+1)\nu)$ one has $y+\dot\psi^{\sharp,1}-\dot\psi^{\sharp,2}=y+r/(2\cc)$ and
we assumed that $\nu\leq\trian/|r/(2\cc)+y|$, it follows that
\begin{align}\notag%\label{eq2150}
\ph^n_1((j^n_1+1)\nu)&=\ph^n_1(\tau^n_1)+(y+r/(2\cc))((j^n_1+1)\nu-\tau^n_1)-(\check{R}^n((j^n_1+1)\nu)-\check{R}^n(\tau^n_1))\\\notag
&>\Barr -\trian -(\check{R}^n((j^n_1+1)\nu)-\check{R}^n(\tau^n_1)).
\end{align}
We now define a new process that starts at time $(j^n_1+1)\nu$, having the form
of the game dynamics with the initial state $\Barr -\trian$, the path
$\psi^*:=\psi^*_{\Barr-\trian}$, and rejection process $\vr^n$.
For $t\ge(j^n_1+1)\nu$ set
\begin{align}\notag%\label{eq2151}
\ph^n_2(t) &=
\Barr -\trian + y(t-(j^n_1+1)\nu) +\psi^{*,1}(t-(j^n_1+1)\nu)\\
 &\notag \qquad - \psi^{*,2}(t-(j^n_1+1)\nu)- \vr^n(t-(j^n_1+1)\nu) ,
\end{align}
where
$\vr^n (s)= \check{R}^n(s )-\check{R}^n(\tau^n_1)$, $s\ge 0$. Notice that at time $t=(j^n_1+1)\nu$ there is an initial amount of rejections
$\check{R}^n((j^n_1+1)\nu )-\check{R}^n(\tau^n_1))\ge 0$, and therefore,
\begin{align}\notag%\label{eq2152}
\ph^n_1((j^n_1+1)\nu)> \Barr -\trian -(\check{R}^n((j^n_1+1)\nu)-\check{R}^n(\tau^n_1))=\ph^n_2((j^n_1+1)\nu).
\end{align}
Therefore
\begin{equation}
  \label{07}
  \ph^n_3:=\by + \psi^{\hs,1} - \psi^{\hs,2} - \check{R}^n\ge\ph^n_2 \quad
  \text{ on the interval } [(j^n_1+1)\nu,\T].
\end{equation}
Let
\begin{equation}\notag%\label{eq2153}
\tau^n_2:= \inf\{t\ge (j^n_1+1)\nu : \ph^n_2(t)\le 0\}.
\end{equation}
From \eqref{eq2047c} it follows that it takes $\int_0^{\Barr-\trian} 1/\sqrt{y^2-h(\xi)/\cc}\;d\xi$ time units for the path
$\Barr-\trian +y\cdot + \psi^{*,1}(\cdot) - \psi^{*,2}(\cdot)$ to reach the level zero. Therefore,
\begin{align}\label{eq2154}
\tau^n_2< T_2\qquad\text{with probability } 1.
\end{align}
Define
\begin{equation}\notag%\label{eq2155}
\ph^n(t)=\left\{\begin{array}{ll}
               \ph^n_1(t)  &\qquad\;\;\quad\ 0\le t \le (j^n_1+1)\nu,
\\
                \ph^n_2(t) & \  (j^n_1+1)\nu < t\leq (j^n_2+1)\nu.
              \end{array}
\right.
\end{equation}

\skp\noi{\bf Step 4 (Estimate on the time $\tau^n_2$):}
Consider the $\calP^{2I}$-path $\upsi^\hs$, which is defined in a similar way to $\upsi^\sharp$.
From Lemma \ref{lem_underline_psi} it follows that there exists $\upsi^\hs=(\upsi^{\hs,1},\upsi^{\hs,2})\in\calP^{2I} $ such that
\begin{align}\label{eq2155b}
(\utheta\cdot\upsi^{\hs,1},\utheta\cdot(\upsi^{\hs,2}\circ\brho)) = (\psi^{\hs,1},\psi^{\hs,2})
\end{align}
and
\begin{align}\label{eq2155c}
\Jr_k(T,\upsi^{\hs,k})= \Ir_k(T,\psi^{\hs,k}) ,\quad k=1,2.
\end{align}
Let
\begin{align}\notag%\label{eq2158}
H^n_k=\{\tau^n_k\in N_{j^n_k} \},\quad k=1,2,
\end{align}
and similarly to Step 2, set
\begin{align}\notag%\label{eq2157}
j^n_2=\underset{j\in\{j^n_1+1,\ldots,j^n_1+\lfloor T_2/\nu\rfloor+2\}} {\arg\max}\PP\left(
\tau^n_2\in N_j \mid \Om^n_{\delta_3,{(j^n_1+1)\nu+T_2}}(\upsi^\hs)  \cap H^n_1\right).
\end{align}
\begin{lemma}\label{lem_n2}
One has
\begin{align}\notag%\label{eq2159}
\liminf_{n\to\iy}\frac{1}{b_n^2}\log\PP\left( H^n_2
\mid \Om^n_{\delta_3,{(j^n_2+1)\nu}}(\upsi^\hs)\cap H^n_1  \right)
\ge -2\eps_1.
\end{align}
\end{lemma}

\skp\noi{\bf Proof:} Recall that $\tau^n_2 \in [(j^n_1+1)\nu,(j^n_1+1)\nu+T_2)$ with probability $1$. Therefore
\begin{align}\notag%\label{eq2160}
\PP\left( H^n_2
\mid \Om^n_{\delta_3,{(j^n_1+1)\nu+\T_2}}(\psi^\hs)\cap H^n_1  \right)\ge
\frac{1}{\left(\lfloor T_2/\nu\rfloor +1\right)}.
\end{align}
The rest of the proof is similar to the proof of Lemma \ref{lem_n1} and is therefore omitted.
\hfill $\Box$

\skp\noi{\bf Step 5 (Bounding the cost from below):}
Let us denote
\[
\Om^n_2=\Om^n_{\delta_3,{(j^n_2+1)\nu}}(\upsi^\hs).
\]
Consider the event $\Om^n_3:=H^n_1\cap H^n_2\cap\Om^n_2\cap E^n$.
On this event, we bound from below the sum
\begin{align}\notag%\label{eq2161}
&\int_0^{\tau^n_2} h(\check X^n(s))ds+r\check R^n(\tau^n_2).%\\\notag
\end{align}
By \eqref{08}, $\check X^n\ge\ph^n_1-2\del_1$ on $[0,(j^n_1+1)\nu]$, and
$\check X^n\ge\ph^n_3-2\del_1$ on $[(j^n_1+1)\nu,\tau^n_2]$
(thanks to the fact that
in Lemma \ref{lem_n1}, $\psi$ is arbitrary on the latter time interval).
Therefore, by \eqref{eq2144} and the definition of $\ph^n$,
on the time interval $[0,\tau^n_2)$,
\begin{align}\notag%\label{eq2162}
\check{X}^n\ge \ph^n-2\delta_1.
\end{align}
Since we chose $\delta_1$ such that $\osc_D(2\delta_1,h)<\eps_1$ it follows that
\begin{align}\label{eq2163}
&\int_0^{\tau^n_2} h(\check{X}^n(t))dt+r\check{R}^n(\tau^n_2)\\\notag
&\quad\ge
\int_0^{\tau^n_2} h(\ph^n(t))dt+r\check{R}^n(\tau^n_2) -\eps_1\tau^n_2\\\notag
&\quad\ge
\Big[\int_0^{\tau^n_1} h(\ph^n_1(t))dt
+r\check{R}^n(\tau^n_1)\Big] + \Big[\int_{(j^n_1+1)\nu}^{\tau^n_2} h(\ph^n_2(t))dt +r(\check{R}^n(\tau^n_2)-\check{R}^n(\tau^n_1))\Big]  -\eps_1\tau^n_2,
\end{align}
where the last inequality follows since $h$ is nonnegative and $\tau^n_1<(j^n_1+1)\nu$.
We now bound from below the three terms above.
From inequality \eqref{eq2045}, the inequality $\tau^n_1\ge j^n_1\nu$, the definitions of $\psi^\hs$ and $\psi^\sharp$, and \eqref{eq2118} it follows that
\begin{align}\label{eq2164}
&\int_0^{\tau^n_1} h(\ph^n_1(t))dt +r\check{R}^n(\tau^n_1)\ge r(x-(\Barr+\trian)) + \Ir(\tau^n_1,\psi^\sharp)\ge r(x-(\Barr+\trian)) + \Ir(j^n_1\nu,\psi^\sharp)
\\\notag
&\quad= r(x-(\Barr+\trian)) + \Ir((j^n_1+1)\nu,\psi^\hs) - \int_{j^n_1\nu}^{(j^n_1+1)\nu}[\cco(\dot\psi^{\sharp,1})^2(t) + \cct(\dot\psi^{\sharp,2})^2(t)]dt\\\notag
&\quad= r(x-(\Barr+\trian)) + \Ir((j^n_1+1)\nu,\psi^\hs) -r^2\nu/(4\cc)\\\notag
&\quad\ge  r(x-(\Barr+\trian)) + \Ir((j^n_1+1)\nu,\psi^\hs) -\eps_1.
\end{align}
To bound the second term notice that
\begin{align}\notag%\label{eq2166}
\bar\ph^n_2(s):=\ph^n_2(s+(j^n_1+1)\nu) = \Barr -\trian + ys +\psi^{*,1}(s) - \psi^{*,2}(s) - \vr^n(s) ,\quad s\ge0.
\end{align}
Denote $\check{\tau}^n_2 = \tau^n_2-(j^n_1+1)\nu$. This is the first time when
$\bar\ph^n_2$ hits zero.
Therefore, from the definitions of $\Ir$, $\check{\tau}^n_2$, and $\psi^\hs$, the inequality $\tau^n_2\ge j^n_2\nu$, and from inequality \eqref{eq2047b} it follows that
\begin{align}\label{eq2167}
&\int_{(j^n_1+1)\nu}^{\tau^n_2} h(\ph^n_2(t))dt +r(\check{R}^n(\tau^n_2)-\check{R}^n(\tau^n_1))=\int_{0}^{\check{\tau}^n_2 } h(\bar\ph^n_2(u))du +r\vr^n(\tau^n_2)\\\notag
&\quad= c(\Barr-\trian,\check{\tau}^n_2,\psi^*, \vr^n)+\Ir(\check{\tau}^n_2,\psi^*) \\\notag
&\quad=  c(\Barr-\trian,\check{\tau}^n_2,\psi^*, \vr^n)+\Ir(\tau^n_2,\psi^\hs) -\Ir((j^n_1+1)\nu,\psi^\hs) \\\notag
&\quad\ge c(\Barr-\trian,\check{\tau}^n_2,\psi^*, \vr^n)+\Ir(j^n_2\nu,\psi^\hs) -\Ir((j^n_1+1)\nu,\psi^\hs)\\\notag
&\quad\ge V(\Barr-\trian)+\Ir(j^n_2\nu,\psi^\hs) -\Ir((j^n_1+1)\nu,\psi^\hs).
\end{align}
Inequality \eqref{eq2154} gives a lower bound on the third term.
Combining inequalities \eqref{eq2163}, \eqref{eq2164}, and \eqref{eq2167} we obtain
\begin{align}\label{eq2168}
&\int_0^{\tau^n_2} h(\check{X}^n(t))dt+r\check{R}^n(\tau^n_2)\\\notag
&\quad\ge
r(x-(\Barr+\trian))+V(\Barr-\trian)-\eps_1(\T+1) +\Ir(j^n_2\nu,\psi^\hs) .
\end{align}

\skp\noi{\bf Step 6 (Bounding the limit from below):} We are now ready to
prove \eqref{eq2107f}.
First, notice that there are only finitely many possible pairs
$\{(j_1,j_2)\in\N^2: 0\le j_1\le j_2\le \T/\nu-1\} $. For each such pair define
\begin{align}\notag%\label{eq2169}
\N_{(j_1,j_2)}=\{n\in\N:(j^n_1,j^n_2)=(j_1,j_2)\}.
\end{align}
If we show that for each pair $(j_1,j_2)$
\begin{equation}\label{09}
\underset{\N_{(j_1,j_2)}}{\liminf}\frac{1}{b_n^2}\log\EE\left[e^{b_n^2\left(\int_0^{\T} h(\check{X}^n(u))du+r\check{R}^n(\T)\right)} \right]\ge V(x)+C_0\eps_1,
\end{equation}
where $C_0$ is a constant independent of $n$ and $\eps_1$, then
\eqref{eq2107f} will follow
on applying Lemma \ref{prop_lower1} and taking $\eps_1\to0$.
Thus in the rest of the proof we focus on a fixed $(j_1,j_2)$, and prove \eqref{09}.
Hereafter, $\liminf$ denotes the limit inferior along the subset.
Denote $\hat\Om^n_2=\Om^n_{\delta_3,(j_2+1)\nu}(\upsi^\hs)$.
From \eqref{eq2168} and since $\tau^n_2\le\T$ it follows that
\begin{align}\label{eq2173}
&\liminf\frac{1}{b_n^2}\log\EE\left[e^{b_n^2\left(\int_0^{\T} h(\check X^n(u))du+r\check R^n(\T)\right)}1_{\{H^n_1\cap H^n_2\cap \hat\Om^n_2\cap E^n \}} \right]\\\notag
&\quad\ge  r(x-(\Barr+\trian))+V(\Barr-\trian) - \eps_1(\T+1) +\Ir(j_2\nu,\psi^\hs) \\\notag
&\qquad +
\liminf\frac{1}{b_n^2}\log\PP\left(H^n_1\cap H^n_2\cap \hat\Om^n_2\cap E^n   \right).
\end{align}
We now estimate the last term above.
By Lemmas \ref{lem_n1} and \ref{lem_n2}, for all $n$ sufficiently large,
\[
\PP\left( H^n_1\mid  \hat\Om^n_2   \right)
\PP\left( H^n_2\mid  H^n_1\cap \hat\Om^n_2   \right)
\ge e^{-5\eps_1 b_n^2}.
\]
Hence
\begin{align}\label{eq2174}
&\liminf\frac{1}{b_n^2}\log\PP\left(H^n_1\cap H^n_2\cap \hat\Om^n_2\cap E^n  \right)\\\notag
&\quad\ge \liminf\frac{1}{b_n^2}\log\left[\PP\left(H^n_1\cap H^n_2\cap \hat\Om^n_2  \right)-\PP\left((E^n)^c\right)\right]\\\notag
&\quad= \liminf\frac{1}{b_n^2}\log\Big[\PP\left( H^n_2\mid  H^n_1\cap \hat\Om^n_2   \right)\\
&\notag
\hspace{9em}\times\PP\left( H^n_1\mid  \hat\Om^n_2   \right)\PP\left(  \hat\Om^n_2   \right)
-\PP\left((E^n)^c\right)\Big]\\\notag
&\quad\ge \liminf\frac{1}{b_n^2}\log\Big( e^{-5\eps_1b_n^2}\,
\PP\left(\hat\Om^n_2\right)
-\PP\left((E^n)^c\right)\Big)\\\notag
&\quad\ge \liminf\frac{1}{b_n^2}\log\Big(
e^{-b_n^2[\Jr((j_2+1)\nu,\bar\psi^\hs)+6\eps_1]} -
e^{b_n^2(-m+1)}\Big)
\\\notag
&\quad= -\Ir((j_2+1)\nu),\psi^\hs) -6\eps_1.
\end{align}
Above, the third inequality follows by Proposition \ref{moderate}
and \eqref{06}.
The last equality uses \eqref{eq2155c} and $m-1>\Ir((j_2+1)\nu),\psi^\hs) +6\eps_1$.
Substituting \eqref{eq2174} in \eqref{eq2173} yields
\begin{align}\label{eq2175}
&\liminf\frac{1}{b_n^2}\log\EE\left[e^{b_n^2\left(\int_0^{\T} h(\check X^n(u))du+r\check R^n(\T)\right)}1_{\{H^n_2\cap H^n_1\cap \tilde\A_{\delta_3,{(j_2+1)\nu}}(\psi^\hs)\cap E^n \}} \right]\\\notag
&\quad\ge  r(x-(\Barr+\trian))+V(\Barr-\trian) -\eps_1(\T+7)+\Ir(j_2\nu,\psi^\hs) -\Ir((j_2+1)\nu),\psi^\hs).
\end{align}

Using \eqref{13} and then \eqref{14} and Proposition \ref{prop_game_value}(iii) gives
\begin{align}\label{eq2176}
&\Ir(j_2\nu,\psi^\hs) -\Ir((j_2+1)\nu),\psi^\hs) = -\int_{j_2\nu}^{(j_2+1)\nu}\left[\cco(\dot\psi^{*,1}_{\Barr-\trian})^2(t)+\cct(\dot\psi^{*,2}_{\Barr-\trian})^2(t)\right]dt\\\notag
&\quad=-\int_{j_2\nu}^{(j_2+1)\nu}\left[\cc(\dot\om^*_{\Barr-\trian})^2(t)\right]dt\\\notag
&\quad=-\int_{j_2\nu}^{(j_2+1)\nu}\left[\cc\left(-y-\sqrt{y^2-h(x+yt+\om^*_{\Barr-\trian}(t))/\cc}\right)^2(t)\right]dt.\\\notag
&\quad\ge -\nu\cc y^2\ge -\eps_1,
\end{align}
where the last two inequalities follow from the negativity of $y$ and \eqref{eq2118}.
From \eqref{eq2175} and \eqref{eq2176} and by recalling that for $x>\Barr$ one has
$V(x)=r(x-\Barr)+V(\Barr)$, it follows that
\begin{align}\notag%\label{eq2177}
&\liminf\frac{1}{b_n^2}\log\EE\left[e^{b_n^2\left(\int_0^{\T} h(\check X^n(u))du
+r\check R^n(\T)\right)}1_{\{H^n_1\cap H^n_2\cap \Om^n_{\delta_3,{(j^n_2+1)\nu}}(\psi^\hs)\cap E^n \}} \right]\\\notag
&\quad\ge  V(x)-\eps_1(\T+8).
\end{align}
This proves \eqref{09}. Hence the result is proved for the case $x>\beta_0$.

Finally, consider $x\le\beta_0$.
The considerations here are simpler than in the previous case.
In case that the initial state is exactly $\Barr$, the decision maker can reject a (small)
amount of $\trian$ at time $t=0$. Then the proof that $V(x)$
is a lower bound requires the focusing only on data near $\psi^*_{\Barr-\trian}$,
starting at time zero. In case that the initial state is lower than $\Barr$,
one uses the same arguments, with data near $\psi^*_x$.
\qed

\subsection{Proof of Lemma \ref{lem_n1}:}
The proof has two parts. On the first we show that for sufficiently large $n$ one has
\begin{align}\label{eq2139}
\PP\left( \tau^n_1\in[0,T_1) \mid \Om^n_{\delta_3,{T_1}}(\psi)  \right)\ge  1/2.
\end{align}
Since the interval $[0,T_1]$ is divided into at most $\lfloor T_1/\nu\rfloor+1$ subintervals,
there exists an interval $N_j$ such that the conditional probability
of $\tau^n_1\in N_j$ is at least $\frac{1}{2\left(\lfloor T_1/\nu\rfloor +1\right)}$.
Thus, as a result of \eqref{eq2139},
\begin{align}\label{15}
\PP\left( \tau^n_1\in N_{j^n_1} \mid \Om^n_{\delta_3,{T_1}}(\psi)  \right)\ge \frac{1}{2\left(\lfloor T_1/\nu\rfloor +1\right)}.
\end{align}
On the second part we use this to deduce \eqref{eq2138}.

\skp

\noi{\bf Part a:} Set
\begin{align}\notag%\label{eq2141}
E^n_1=\{\tau^n_1\ge T_1\}.
\end{align}
Write $\Om^n_{\delta_3,{T_1}}(\psi)$ as $\Om^n$.
We analyze the event $\Om^n_1:=\Om^n\cap E^n\cap E^n_1$.
On this event,
we bound from below the process $\check{X}^n$ on the time interval $[0,T_1]$ and the total number of rejections until time $T_1$.
By the triangle inequality it follows that for every $i\in\calI$ and every $n$
\begin{equation}\label{eq2142}
\|\tilde S^n_i\circ T^n_i- \psi^2_i\circ\brho_i\|_{\T}
\leq \|\tilde S^n_i\circ T^n_i
- \psi^2_i\circ T^n_i\|_{\T}+\|\psi^2_i\circ T^n-\psi^2_i\circ\brho_i\|_{\T}\le 2\del_2,
\end{equation}
where we have bounded each of the terms on the r.h.s.\ by $\delta_2$;
the bound of the first term follows by \eqref{eq2130},
and the bound of the second
follows by the definition of $E^n$ and from \eqref{eq2129}. Similarly,
\begin{equation}\label{eq2142b}
\|\tilde A^n_i-\psi^1_i\|_{\T}\leq \delta_2,\quad i\in\calI.
\end{equation}
Since $\utheta^n\to\utheta$ it follows from \eqref{eq2128b}, \eqref{eq2142}, and \eqref{eq2142b} that for sufficiently large $n$,
\begin{align}\label{eq2142c}
\|\check{\bS}^n \circ T^n-\psi^2\circ\brho\|_{\T}&=\|\utheta^n\cdot\tilde\bS^n\circ T^n
- \utheta\cdot\psi^2\circ\brho\|_{\T}<\frac{\delta_1}{2},\\\label{eq2142d}
\|\check{A}-\psi^1\|_{\T}&=\|\utheta^n\cdot\tilde{A}^n
-\utheta\cdot\psi^1\|_{\T}<\frac{\delta_1}{4}.
\end{align}
Moreover, for sufficiently large $n$ one has, for $t\in[0,\bar T]$,
$|\check{X}^n(0)+\check{y}^nt-\by(t)|\leq\delta_1/4$. Using the above inequalities it follows that for every $u\in[0,T_1]$ one has
\begin{align}\label{08}
\check{X}^n(u)
&=\check{X}^n(0) + \check{y}^nt +\check{A}^n(u)-\check{\bS}^n\circ T^n(u)+ \check{Z}^n(u) - \check{R}^n(u) \notag \\
&\ge \ph^n_1(u)-\delta_1 - \sup_{t\in[0,\bar T]}\|\check{X}^n(0)+\tilde y^nt-\by(t)\|
- \|\check{A}^n - \psi^1\|_{\T}
- \|\check{\bS}^n\circ T^n - \psi^2\|_{\T}\notag \\
&\ge \ph^n_1(u)-2\delta_1.
\end{align}
By the definition of $E^n_1$
and the choice of $\delta_1$ (see \eqref{eq2126}),
\begin{equation}\label{eq2143}
\check{X}^n(u)\ge \Barr+\trian-2\delta_1\ge \Barr+\trian/2.
\end{equation}
By using similar arguments and the inequality $\check{X}^n\le D^n$ one obtains
\begin{align}\label{eq2144}
\check{R}^n(T_1)&\ge \by(T_1)-D^n -\delta_1
+ \psi^1(T_1)-\psi^2(T_1) = x-D^n -\delta_1+(y+r/(2\cc))T_1\\\notag
&\ge x-D^n -1+(y+r/(2\cc))T_1\ge x-D -2+(y+r/(2\cc))T_1 ,
\end{align}
where the equality follows by the definition of $\psi^\sharp$, the second inequality follows by the choice of $\delta_1$, and the last inequality follows since $\limn D^n=D$.
From \eqref{eq2107a}, \eqref{eq2119}, and since $T_1<\T$,
we obtain\footnote{Recall the convention $\eps=0$ and generic $r$.}
for sufficiently large $n$
\begin{align}\notag%\label{eq2144b}
V(x)+1>\frac{1}{b_n^2}\log\EE\left[e^{b_n^2\left(\int_0^{T_1} h(\check{X}^n(u))du+r\check{R}^n(T_1)\right)} \right].
\end{align}
Along with \eqref{eq2143} and \eqref{eq2144} it follows that
\begin{align}\notag%\label{eq2145}
V(x)+1&>\frac{1}{b_n^2}\log\EE\left[e^{b_n^2\left(\int_0^{T_1} h(\check{X}^n(u))du+r\check{R}^n(T_1)\right)}1_{\Om^n_1} \right]\\\notag
&\ge
\frac{1}{b_n^2}\log\EE\left[e^{b_n^2\left((h(\Barr+\trian/2)+yr+r^2/(2\cc))T_1
+r(x-D-2)\right)}1_{\Om^n_1} \right]\\\notag
&=r(x-D-2) +(h(\Barr+\trian/2)-h(\Barr) +r^2/(4\cc))T_1 \\
& \notag
\qquad + \frac{1}{b_n^2}\log\PP(\Om^n_1).
\end{align}
The above equality follows since
 $r^2/(4\cc)+ry+h(\Barr)=0$, which in turn follows since $x>\Barr$ and therefore $\Barr<D$ .
Since
\begin{align}\notag%\label{eq2146}
\PP(\Om^n_1)=
\PP(E^n_1\cap\Om^n\cap E^n)
\ge \PP(E^n_1 | \Om^n)\PP(\Om^n)-\PP(( E^n)^c),
\end{align}
and using \eqref{06}, it follows that
\begin{align}\notag%\label{eq2147}
\PP(E^n_1 | \Om^n) &\le e^{b_n^2 \left(V(x)+1+r(2+D-x)-T_1(h(\Barr+\trian/2)-h(\Barr)-r^2/(4\cc))\right)}
\,\PP(\Om^n)^{-1}
\\ \notag
&\quad+ e^{-m b_n^2}\,\PP(\Om^n)^{-1}.
\end{align}
We show that for sufficiently large $n$,
each of the terms on the r.h.s.\ can be bounded by $1/4$. From
Proposition \ref{moderate} and \eqref{eq2128c}, it follows that, for all large $n$,
\begin{align}\notag%\label{eq2148}
\frac{1}{b_n^2}\log\PP(\Om^n)&\ge -\inf_{\upsi\in\A_{\delta_3,{T_1}}(\psi)}\Jr(T_1,\upsi) -1
\ge  -\Jr(T_1,\psi) -1 \\ \notag
&= -\Ir(T_1,\psi^\sharp) -1= -r^2T_1/(4\cc)-1 .
\end{align}
Hence by the definition of $T_1$, the first term is bounded by $1/4$.
Since $m>\Ir(\T,\psi^\sharp)$ and $\T>T_1$,
so is the second term. As a result, \eqref{eq2139} holds.

\skp

\noi{\bf Part b:}
By the definition of an admissible control, $R^n$ is adapted to the filtration
$\calF_t:=\sig\{A^n_i(u),S^n_i(T^n_i(u)),i\in\calI,u\le t\}$, hence so is $\check R^n$,
and, by \eqref{101}, so is $\ph^n_1$. Since $T^n_i(t)\le t$, $t\ge0$, and $T^n_i$ are
themselves adapted to $\calF_t$, it follows that
the event $\{\tau^n_1\in N_{j^n_1}\}$ is measurable on $F^n_{\nu(j^n_1+1)}$,
where $F^n_t=\sig\{A^n_i(u), S^n_i(u),i\in\calI,u\le t\}$. Note that
$F^n_t=\sig\{\tilde A^n_i(u), \tilde S^n_i(u),i\in\calI,u\le t\}$.

Fix $v>0$ and a sequence $v_n$, with $v_n< v$ (both deterministic).
We will show the following.
Given a constant $c_1>0$ and a sequence of events $Q_n\in F^n_{v_n}$, for every
$\eps>0$ there exists $\del>0$ and $n_1\in\N$, such that
\begin{equation}
  \label{18}
  p_1^n:=\PP(Q_n|\Om^n_{\del,v}(\psi))\ge c_1, \quad n\ge1,
\end{equation}
implies
\begin{equation}
  \label{19}
  p_2^n:=\PP(Q_n|\Om^n_{\del,v_n}(\psi))\ge e^{-\eps b_n^2},\quad n\ge n_1.
\end{equation}
Note that this will prove \eqref{eq2138}, based on \eqref{15} that has now been established
by part (a).

Extending the definition of $\Om^n_{\del,t}(\psi)$
\eqref{eq2134}, we let
\[
\Om^n_{\del,a,b}(\psi)=
\{\sup_{s\in[a,b]}\|(\tilde A^n,\tilde S^n)(u)-\psi(u)\|<\del\},
\]
for $0\le a\le b$.
Also, we drop $\psi$ from the notation $\Om^n_{\del,a}$ and $\Om^n_{\del,a,b}$.
Note that there is no loss of generality in proving the statement for
$\tilde A^n$ (a collection of $I$ independent renewal processes)
in place of $(\tilde A^n, \tilde S^n)$ (a collection of
$2I$ such processes). Thus we will consider only the former.

To prove the aforementioned statement, let $\eps>0$ be given.
Consider the quantities $p^n_1$ and $p^n_2$, depending on $\del$. Assume that \eqref{18}
is valid. Then we can write
\[
p_1^n=\frac{\PP\{Q_n\cap\Om^n_{\del,v}\}}{\PP\{\Om^n_{\del,v}\}}
=\frac{\PP\{Q_n\cap\Om^n_{\del,v_n}\cap\Om^n_{\del,v_n,v}\}}
{\PP\{\Om^n_{\del,v_n}\cap\Om^n_{\del,v_n,v}\}}.
\]
A basic independence property for a renewal process, to be used, is the following.
Let $A$ be a renewal process of the form
\[
A(t)=\sup\Big\{l\ge0:\sum_{k=1}^l U(k)\le t\Big\},\quad t\ge0,
\]
where $\{U(k)\}$ are iid (compare with \eqref{20}). Fix $t$ and let $\pi$
denote the time of the first jump at or after $t$, namely $\pi=\inf\{u\ge t:A(u)>A(t-)\}$
(convention: $A(0-)=0$).
Consider an event $Q$ measurable on $\sig\{A(u):0\le u\le t\}$. Then, for each $k\in\Z_+$,
the event $Q\cap\{A_t=k\}$ is statistically
independent of the sequence $\{U_{k+1},U_{k+2},\ldots\}$. Based on this
it is not hard to see that, if we let $SA$ denote the shifted version
$SA(u)=A(\pi+u)-A(\pi)$, $u\ge0$, of $A$, we have that $Q$ is independent
of $SA$. For a collection of independent renewal processes, a similar
statement holds if each of them is shifted according to its own first jump after $t$.
To state this property for the processes $(A^n_i)$, if $Q$ is measurable
on $F^n_t$, then it is independent of $(SA^n_i)$. Now, let us apply this
to study $(\tilde A^n)$. Let $\pi_n=(\pi^{n}_i)$
be defined by $\pi^{n}_i=\inf\{u\ge v_n:\tilde A^n_i(u)>\tilde A^n_i(v_n-)\}$. If
$\|\pi_n\|:=\max_i|\pi^n_i-v_n|$, then, given any $k>0$,
$\PP\{\|\pi_n\|>\del\}\le e^{-kb_n^2}$ for all sufficiently large $n$,
as can be verified using the exponential moment assumption and applying Chebychev's inequality.
Hence, given any $k$, for all large $n$,
\[
p^n_{1,1}:=\PP\{Q_n\cap\Om^n_{\del,v_n}\cap\Om^n_{\del,v_n,v}\}\le
\PP\{Q_n\cap\Om^n_{\del,v_n}\cap\Om^n_{\del,v_n,v}\cap\{\|\pi_n\|\le\del\}\}+e^{-kb_n^2}.
\]
Let us denote by $\om$ a modulus (by which we mean a function mapping $\R_+$ to itself
with $\om(0+)=0$), that dominates the modulus of continuity of $\bar\psi^m$
for all $m\le M$.
Then, using the definition of $\Om^n_{\del,v_n,v}$, adding and subtracting
the shifted version $S\tilde A^n$ and using the triangle inequality gives
\[
p^n_{1,1}\le
\PP\{Q_n\cap\Om^n_{\del,v_n}
\cap\{
\|S\tilde A^n-(\psi(v_n+\cdot)-\psi(v_n))\|_{v-v_n-\del}<\del'\}\}+e^{-kb_n^2},
\]
where $\del'=2\del+2\om(\del)$.
Using the independence alluded to above,
\[
p^n_{1,1}\le
\PP\{Q_n\cap\Om^n_{\del,v_n}\}
\PP\{\|S\tilde A^n-(\psi(v_n+\cdot)-\psi(v_n))\|_{v-v_n-\del}<\del'\}
+e^{-kb_n^2}
\]
Shifting back gives
\[
p^n_{1,1}\le\PP\{Q_n\cap\Om^n_{\del,v_n}\}
\PP\{\Om^n_{\del'',v_n,v-\del}\}+e^{-kb_n^2}
\]
for $\del''$ that can be made arbitrarily small by taking $\del$ to be small.
A similar argument shows that, for all large $n$,
\[
p^n_{1,2}:=\PP\{\Om^n_{\del,v_n}\cap\Om^n_{\del,v_n,v}\}\ge
\PP\{\Om^n_{\del''',v_n}\}\PP\{\Om^n_{\del''',v_n,v}\}-e^{-kb_n^2},
\]
for suitably chosen $\del'''>0$, that again, can be made arbitrarily small by taking small $\del$.
Thus
\begin{align*}
p^n_1&=\frac{p^n_{1,1}}{p^n_{1,2}}
\le
\frac{\PP\{Q_n\cap\Om^n_{\del,v_n}\}\PP\{\Om^n_{\del'',v_n,v-\del}\}+e^{-kb_n^2}}
{\PP\{\Om^n_{\del''',v_n}\}\PP\{\Om^n_{\del''',v_n,v}\}-e^{-kb_n^2}} \\
&=
\frac{\PP\{Q_n|\Om^n_{\del,v_n}\}\PP\{\Om^n_{\del,v_n}\}\PP\{\Om^n_{\del'',v_n,v-\del}\}+e^{-kb_n^2}}
{\PP\{\Om^n_{\del''',v_n}\}\PP\{\Om^n_{\del''',v_n,v}\}-e^{-kb_n^2}}.
\end{align*}
Using the LDP and writing $\Ir[a,b]$ for $\sum_i\int_a^b\cc_i\dot\psi(u)^2du$ (see \eqref{eq2018})
gives
\begin{align*}
p^n_1&\le
\frac{p^n_2e^{b_n^2(-\Ir[0,v_n]+\eps')}e^{b_n^2(-\Ir[v_n,v-\del]+\eps')}+e^{-kb_n^2}}
{e^{b_n^2(-\Ir[0,v_n]-\eps')}e^{b_n^2(-\Ir[v_n,v]-\eps')}-e^{-kb_n^2}} \\
&=
\frac{p^n_2e^{b_n^2(-\Ir[0,v-\del]+2\eps')}+e^{-kb_n^2}}
{e^{b_n^2(-\Ir[0,v]-2\eps')}-e^{-kb_n^2}},
\end{align*}
for arbitrary $\eps'>0$, provided that $\del$ is sufficiently small.
Hence by selecting $k$ sufficiently large,
\begin{align*}
p^n_1&\le
\frac{p^n_2e^{b_n^2(-\Ir[0,v-\del]+2\eps')}+e^{-kb_n^2}}
{\frac{1}{2}e^{b_n^2(-\Ir[0,v]-2\eps')}} \\
&=
2p^n_2e^{b_n^2(4\eps'+\Ir[v-\del,v])}+2e^{-kb_n^2}e^{b_n^2(\Ir[0,v]+2\eps')}.
\end{align*}
Hence, again by selecting $k$ large, for all large $n$,
\[
p^n_2\ge \frac{1}{2}c_1e^{-b_n^2(4\eps'+\Ir[v-\del,v])}-e^{b_n^2(-k+\Ir[0,v-\del]-2\eps')}
\ge\frac{1}{4}c_1e^{-b_n^2(4\eps'+\Ir[v-\del,v])}.
\]
Selecting $\del>0$ such that $4\eps'+\Ir[v-\del,v]<\eps$ gives \eqref{19} for some $n_1$.
This completes the proof of part (b), and the lemma.
\qed

\section{A nearly optimal policy}\label{sec_upper}
\beginsec

In this section we show that the policy from \cite{ata-shi}
is AO for the present setting. While the policy is similar,
the proof of AO is quite different, as the paper \cite{ata-shi} addresses the diffusion scale,
rather than the MD scale.

Let the classes be labeled so that
$h_1\mu_1\ge h_2\mu_2\ge\cdots\ge h_I\mu_I$.
Let $\gamma:[0,D]\to\calX$ be a Borel measurable mapping satisfying
\begin{equation}\label{103}
\gamma(w)\in\argmin_{\xi}\{\uh\cdot\xi:\xi\in\calX,\utheta\cdot\xi=w\},\qquad w\in[0,D].
\end{equation}
We note on passing that, as shown in \cite[Theorem A.1]{ata-coh},
one can equivalently work with one-dimensional dynamics, thanks to the fact that
the minimum over queue length $\xi$ in the above expression
is a function of only of the (one-dimensional) workload $w$.

Since the mapping $\uxi\mapsto\uh\cdot\uxi$ is linear and the domain $\calX$ is polyhedral,
it can be assumed, without loss of generality, that
$\gamma$ is continuous and takes values on the boundary of $\calX$.
We have, by definition, that $\utheta\cdot\gamma(w)=w$, and
$\uh\cdot\gamma(w)= h(w)\le \uh\cdot\xi$ for
$\xi\in\calX$ for which $\utheta\cdot\xi=w$.
A particular selection of $\gamma$ is as follows.
Given $w\in[0,D]$, set $(j,\xi)=(j,\xi)(w)$ by
$w\in[\hat D_j,\hat D_{j-1})$
and $\xi=\xi(w):=(w-\hat D_j)/\theta_j$,
where
\[
\hat D_j:=\sum_{i=j+1}^I\theta_iD_i,\quad j\in\{0,\ldots,I\}
\]
and one has $0=\hat D_I<\hat D_{I-1}<\cdots<\hat D_1<\hat D_0=\utheta\cdot \uD=D$, $\uD:=(D_1,\ldots,D_I)$. Then
\begin{align}\label{eq2303}
\gamma(w)=\sum_{i=j+1}^ID_i\ue_i+\xi \ue_j.
\end{align}

Approximate $\gamma$ by a curve
that is bounded away from the part of the boundary of $\calX$ that corresponds to
the buffer limit, namely $\pl\calX=\{x\in\calX:x_i=D_i \text{ for some } i\}$.
Fix $0<\eps_0<\min_iD_i/4$.
Let $a_i=D_i-3\eps_0$, $i\in\calI$, and
$a^*=\Barr\w(\utheta\cdot a)<D$. Note that if $\eps_0$ is small
then $a^*=\Barr$ (unless $\Barr=D$).
Define $\gamma^a[0,D]\to\calX$
first on $[0,\utheta\cdot a]$ as the function obtained upon replacing
the parameters $(D_i)$ by $(a_i)$ in \eqref{eq2303}.
That is, for $w\in[0,\utheta\cdot a)$, the variables $j=j(w)$ and $\xi=\xi(w)$ are determined
via
\begin{equation}\label{eq2304}
w=\sum_{i=j+1}^I\theta_ia_i+\theta_j\xi, \qquad j\in\{1,2,\ldots,I\},\,\xi\in[0,a_j),
\qquad
\gamma^a(w)=\sum_{i=j+1}^Ia_i\ue_i+\xi \ue_j.
\end{equation}
Given $w\in[0,\utheta\cdot a)$, we will sometimes refer to the unique pair
$(j,\xi)$ alluded to above as the {\it representation $(j,\xi)$ of $w$ via \eqref{eq2304}}.
Next, on $[\utheta\cdot a,\utheta\cdot \uD]$
define $\gam^a$ as the linear interpolation between the points $(\utheta\cdot a,a)$
and $(\utheta\cdot \uD,\uD)$.
Also set $\hat a_j=\sum_{i=j+1}^I\theta_ia_i$, $j\in\{0,1,\ldots,I\}$. Let
\[
h^a(w):=\min\{\uh\cdot\xi:\xi\in\calX, \utheta\cdot\xi=w,
\xi_i\le\gamma^a_i(w),i=1,\ldots,I\}=\uh\cdot\gamma^a(w),\qquad
w\in[0, \utheta\cdot a].%\label{eq2307}
\]
Note the similarity to the payoff $h$ in \eqref{eq2026g}. Note that the construction
depends on the parameter $\eps_0$, and denote
\begin{align}\label{eq2307b}
\om_1(\eps_0)=\sup_{[0,\utheta\cdot a]}| h^a- h|.
\end{align}
By the choice of $a$ it is clear that $\om_1(0+)=0$.

Before providing the precise construction of the policy,
we explain its rationale.
The solution to the DG indicates that rejections should occur when the normalized workload
in the system is above the threshold $\beta_0$,
and that most rejections should be from a specific class, $i^*$, defined in \eqref{eq20238aa}.
The DG solution also indicates that prioritization should be according to \eqref{103}
(see the proof of (100) in \cite[Theorem A.1]{ata-coh}) and that, consequently,
the resulting normalized queue length processes should be close to the curve $\gamma$.
These two goals are contradictory, as parts of the curve $\gamma$ lie on the
part $\pl\calX$ of the domain where some of the buffers are full,
and so even small stochastic fluctuations
cause rejections due to the buffer size constraints. Such rejections do not satisfy
the requirement to reject only when the workload is above the specified threshold,
nor that rejections are from class $i^*$. To address this issue,
we have defined the curve $\gamma^a$, which approximates $\gamma$
without intersecting the part $\pl\calX$ of the boundary.
The service policy is designed to keep the normalized queue length
processes close to this curve.

The precise definition of the policy is provided by specifying $(B^n(t),R^n(t))$
as a function of $X^n(t)$.

{\it Rejection policy:}
As under any policy, in order to meet the
buffer size constraint \eqref{10}, all \emph{forced rejections} take place. That is,
if a class-$i$ arrival occurs at a time $t$ when $\tilde X^n_i(t-)+\tfrac{1}{b_n \sqrt{n}}>D_i$,
then it is rejected.
Apart from that, no rejections occur from any class except class $i^*$, which is defined through \eqref{eq20238aa},
and no rejections occur (from any class) when $\utheta\cdot \tilde X^n<a^*$.
When $\utheta\cdot \tilde X^n\ge a^*$, all class-$i^*$ arrivals are rejected,
and these rejections are called \emph{overload rejections}.

{\it Service policy:}
For each $\bar x\in\calX$ define the class of low priority
\[
\calL(\bar x)=\max\{i:x_i<a_i\},%\label{eq2308}
\]
provided $x_i<a_i$ for some $i$, and set $\calL(\bar x)=I$ otherwise.
The complement set is the set of high priority classes:
\[
\calH(\bar x)=\calI\setminus\{\calL(\bar x)\}.%\label{eq2309}
\]
When there is at least one class among $\calH(\bar x)$
having at least one customer in the system,
$\calL(\bar x)$ receives no service, and
all classes within $\calH(\bar x)$, having at least one customer,
receive service at a fraction proportional to their
traffic intensities. Namely, denote $\calH^+(\bar x)=\{i\in \calH(\bar x):x_i>0\}$,
and define $\rho'(\bar x)\in\R^I$ as
\begin{equation}\label{eq2310}
\rho'_i(\bar x)=\begin{cases}
0, & \text{if } \bar x=\uzero,\\
\ds\frac{\rho_i1_{\{i\in\calH^+(\bar x)\}}}{\sum_{k\in\calH^+(\bar x)}\rho_k}, &
\text{if } \calH^+(\bar x)\ne\emptyset,\\
e_I,& \text{if $x_i=0$ for all $i<I$ and $x_I>0$.}
\end{cases}
\end{equation}
(Note that $\calH^+(\bar x)=\emptyset$ can only happen if $x_i=0$ for all $i<I$,
which is covered by the first and last cases in the above display).
Then for each $t$,
\begin{equation}\label{eq2311}
  B^n(t)=\rho'(\tilde X^n(t)).
\end{equation}
Note that when $\calH^+(\bar x)\ne\emptyset$,
\begin{equation}\label{eq2312}
  \rho'_i(\bar x)>\rho_i\quad \text{for all } i\in \calH^+(\bar x).
\end{equation}
That is, all prioritized classes receive
a fraction of effort strictly greater than the respective traffic intensity.
Also note that $\sum_iB^n_i=1$ whenever $\tilde X^n$ is nonzero. This is therefore
a work conserving policy.

\begin{theorem}\label{thm_upper}
  Let Assumption \ref{assn1} hold. For every $\eps_0>0$ and $n\in\N$, denote the policy constructed above by $U^n(\eps_0)$. Then, for all sufficiently large $T$,
  \[
  \limsup_{n\to\iy}J^n(T,\tilde X^n(0),U^n(\eps_0))\le V(x)+\omega(\eps_0),
  \]
  where $\om:\R_+\to\R_+$ is a function satisfying $\omega(0+)=0$.
\end{theorem}

\noi{\bf Proof of Theorem \ref{thm_upper}:}
Introduce the notation
\[
H^n_t=\int_0^t\uh\cdot\tilde X^n(u)du+\ur\cdot\tilde R^n(t). %\label{eq2328}
\]
Fix $T>0$ sufficiently large for the identity $V(T,x)=V(x)$ stated in
Proposition \ref{prop_game_value}(i) to hold. First, notice that
\begin{equation}\label{eq2313}
\E\int_0^T e^{b^2_n H^n_t}dt
\leq
T \E\Big[ e^{b^2_nH^n_T}\Big].
\end{equation}
The argument will be based on a bound on the r.h.s.\ of \eqref{eq2313}.

Recall the definition \eqref{eq2018} of $\Jr$, and for ${J}>0$, define
\begin{equation}\notag%\label{eq2314}
\AC_{J}=\{\psi\in \calD([0, T], \R^{2I})\ : \ \Jr(T,\psi)\leq {J}\}.
\end{equation}
Then $\AC_{J}$ is compact in the $J_1$ topology, and
consists of absolutely continuous paths starting at zero.
Fix $0< \eps_1 <\theta_{\min}\eps_0/8$.
Fix also $\del_1>0$ such that
\begin{equation}\label{eq2316b}
\delta_1<\min\{\eps_0/12,\eps_1/(11\ee), \eps_0/(5\ee), \osc_D(\eps_0/(5\ee),h)\},
\end{equation}
where $\ee$ is the constant from Lemma \ref{lem_appendix1}.
For $0<\del<t\le T$ and $\psi\in\calD([0, T], \R^{2I})$, denote
\[
\A_{\delta,t}(\psi)=\{\upsi\in \calD([0,T], \R^{2I}) \ :\ \|\upsi-\psi\|_t<\delta\}
\]
(where we slightly modified the notation \eqref{eq2130}).
By the compactness of $\AC_J$ and the continuity of its members,
one can find a finite number of members $\upsi^1, \upsi^2,\ldots,$ $\upsi^M$
of $\AC_{J}$, and positive constants $\delta^1,\ldots, \delta^M$ with $\delta^m<\delta_1$, satisfying
$\AC_{J}\subset\cup_{m=1}^M \A^m_T$, and
\begin{equation}\label{eq2323}
\inf\{\Jr(T,\upsi):\upsi\in \oo{\A^m_T}\}\geq \Jr(T,\upsi^m)-\eps_0,
\qquad m=1,2,\ldots,M,
\end{equation}
where, throughout, for $0\leq t\leq T$,
$$
\A^m_t=\A_{\del^m,t}(\bar\psi^m).
$$

By the continuity of each of the paths $\bar\psi^m$,
one can find $\nu_1>0$ such that for $m=1,\ldots,M$,
\begin{equation}\label{eq2317}
\osc_T(\nu_1,\bar\psi^m_i)\leq\frac{\delta_1}{2(\theta_{\max}\vee1)\sqrt{I}},
\quad i\in\calI.
\end{equation}

In this proof, $C_{1}, C_{2},\ldots$ denote positive constants that do not depend on
$n,\eps_0,\delta^m$ or ${J}$.
Write $\upsi^m\in\calP^{2I}$ as $(\upsi^{m,1},\upsi^{m,2})$.
For $m=1,\ldots,M$ let
\[
\ph^m= \by+\psi^{m,1}-\psi^{m,2}+\zeta^m-\vr^m,%\label{eq2325}
\]
where, as before, $\by(t)=x+yt$,
\begin{align}\notag%\label{eq2325b}
\psi^{m,1}=\utheta\cdot \upsi^{m,1},
\quad   \psi^{m,2}=\utheta\cdot \upsi^{m,2},
\end{align}
and
\begin{align}\label{eq2326}
(\ph^n,\zeta^m,\vr^m)
=\Gam_{[0,\Barr]}[ \by + \psi^{m,1}- \psi^{m,2}].%\label{eq2327}
\end{align}
Denote $\La^n=\|\tilde{A}^n\|_T+\|\tilde{S}^n\|_T$.
As argued in \cite{ata-shi}, at the bottom of page 595,
\begin{align}\label{eq2329}
\|\tilde R^n (T)\|\leq C_{1}(1+\La^n).
\end{align}
Since $\tilde X^n$ is bounded, one has $H^n_T\leq C_{2}(1+\La^n)$.
Hence, given any ${J}_1>0$,
\[
H^n_T> {J}_1
\quad\text{implies}\quad \La^n > C_{2}^{-1}{J}_1-1=:G({J}_1).%\label{eq2330}
\]
Therefore
\begin{align}\label{eq2331}
\notag
\E[e^{b_n^2H^n_T}] &\le \E[e^{b_n^2[H^n_T\wedge {J}_1]}]
+ \E[e^{b_n^2H^n_T}1_{\{H^n_T>{J}_1\}}]
\\
&\le A^n_1+A^n_2+A^n_3,
\end{align}
where, with $\B=(\cup_{m=1}^M\A^m_T)^c$,
\begin{align*}
&A^n_1=\sum_{m=1}^M\E[e^{b_n^2[H^n_T\wedge {J}_1]}1_{\Om^{n,m}}],
\qquad
\Om^{n,m}=\{(\tilde{A}^n,\tilde{S}^n)\in \A^m_T\},
\\
&
A^n_2=\E[e^{b_n^2[H^n_T\wedge {J}_1]}1_{\{(\tilde{A}^n,\tilde{S}^n)\in \B\}}],\\
&A^n_3=\E[e^{b_n^2C_{2}(1+\La^n)}1_{\{\La^n>G({J}_1)\}}].
\end{align*}
An argument to be presented shortly will show that
there exist $t_1,\ldots,t_M\in[0,T]$ such that for large $n$,
\begin{equation}\label{eq2332}
A^n_1\leq M\max_{m=1}^M
e^{b_n^2[ \int_0^{t^m} h(\ph^m(u))du+r(1+\eps_0)\vr^m(t^m) - \Ir(t^m,\psi^m)
+\om_2(\eps_0)]}+\eps_0,
\end{equation}
where $\om_2(0+)=0$
(this step translates
the multidimensional formulation, by which $H^n_T$ is defined, into
a one-dimensional form, given by $\ph^m$).
As for $A^n_2$ and $A^n_3$, note, by Proposition \ref{moderate}, that for large $n$,
\begin{align}\notag%\label{eq2382}
\frac{1}{b_n^2}\log\p((\tilde{A}^n,\tilde{S}^n)\in \B)\leq -\inf_{\psi\in \B}\Ir(T,\psi)+\eps_0.
\end{align}
Along with the fact that $\B\subset \AC_{J}^c$ and the definition of $\AC_{J}$, this shows
\begin{equation}
\label{eq2333}
A^n_2\leq e^{b_n^2[{J}_1-{J}+\eps_0]}.
\end{equation}
Also, $A^n_3\leq\E[e^{b_n^2((C_{2}+1)\La^n+C_{2}-G({J}_1))}]$.
As shown in the appendix of \cite{ata-bis}, Assumption \ref{assn1} implies that, for any $K<\iy$,
\begin{equation}\notag
\limsup_{n\to\iy}\frac{1}{b^2_n}\log\E[e^{b^2_nK(\|\tilde{A}^n\|_T+\|\tilde{S}^n\|_T)}]
<\infty.%\label{025}
\end{equation}
Hence there exists a constant $C_{3}$ such that
\begin{equation}
\label{eq2334}
A^n_3\le e^{b_n^2 C_{2}+C_{3}-G({J}_1)}.
\end{equation}
Combining \eqref{eq2331}, \eqref{eq2332}, \eqref{eq2333} and \eqref{eq2334},
\begin{align*}\notag%\label{eq2335}
&\limsup\frac{1}{b_n^2}\log\E[e^{b_n^2H^n_T}]
\\\notag%\label{eq202}
& \leq \max_{1\leq m\leq M}\Big[\int_0^{t^m} h(\ph^m(u))du+r(1+\eps_0)\vr^m(t^m)-\Ir(t^m,\psi^m)+\om_2(\eps_0)\Big]\\
&\qquad \vee [{J}_1-{J}+\eps_0]\vee[C_{2}+C_{3}-G({J}_1)]
\\\notag%\label{eq203}
&\leq \sup_{\psi\in \calP,t\in[0,T]}[c^{\eps_0}(x,t,\psi,\al_\Barr[\psi])+\om_2(\eps_0)]\vee [{J}_1-{J}+\eps_0]
\vee[C_{2}+C_{3}-G({J}_1)],\notag%\label{eq204}
\end{align*}
where the cost $c^{\eps_0}$ is defined as $c$ with the rejection cost $r(1+\eps_0)$ instead of $r$.
Now, let $\eps_0\to 0$ first, then ${J}\to\infty$, recalling that
$C_{2}$, $C_{3}$ and $G$ do not depend on
${J}$. Finally let ${J}_1\to\infty$, so $G({J}_1)\to\infty$, to obtain
$$
\limsup V^n(T,\tilde X^n(0))\le
\limsup\frac{1}{b_n^2}\log\E[e^{b_n^2H^n_T}]\leq
\sup_{\psi\in \calP^2,t\in[0,T]} c(x,t,\psi,\al_\Barr[\psi])=V(x),
$$
where in the first inequality we used \eqref{eq2313},
and in the equality we used
the optimality of the $\Barr$-barrier strategy in the game, see Proposition
\ref{prop_game_optimal_strategy}, as well as Proposition \ref{prop_game_value}(i).

We thus turn to the proof of \eqref{eq2332}.
We argue in two steps. In step 1, we show the multidimensional process
$\tilde X^n$ lies close to the minimizing curve. Consequently, we also deduce that no forced
rejections occur, provided $n$ is sufficiently large.
In step 2, we deduce \eqref{eq2332} from step 1.

\skp\noi{\bf Step 1.}
We show that for large $n$,
\begin{equation}\label{eq2338a}
\max_{i}\|\Del^n_i\|_T\le\eps_0,
\end{equation}
where we denote the difference process
\begin{equation}\label{eq2338b}
\Del^n_i(t)=\tilde X^n_i(t)-\gamma^a_i(\check{X}^n(t)),\quad t\in[0,T].
\end{equation}
Denote by $\calG=\{x\in\calX:\utheta\cdot x\le a^*,x=\gamma^a(\utheta\cdot x)\}$
the set of points lying on the minimizing curve, and recall
$\pl^+\calX:=\{x\in\calX:x_i=b_i \text{ for some } i\}$, the set corresponding to the
buffer limit boundary. By the choice of $a$ and $\eps_0$ it follows that $\calG^{\eps_0}$
and $(\pl^+\calX)^{\eps_0}$ do not intersect, where for a set $A\in\R^I$ we denote
$A^{\eps}:=\{x:\dist(x,A)\le\eps\}$.
Forced rejections occur only at times when $\tilde X^n$ lies in $(\pl^+\calX)^{\eps_0}$
(for all $n$ large). As a result, as long as the process
$\tilde X^n$ lies in $\calG^{\eps_0}$, no forced rejections occur.
Thus by showing \eqref{eq2338a}, one also obtains
\begin{align}\label{eq2340}
\tilde R^n(T)=\tilde R^n_{i^*}(T)\ue_{i^*}.
\end{align}
Denote
\begin{align}\notag%\label{eq2341}
\tau^n=\inf\{t\ge 0:\max_{i}|\Del^n_i(t)|\ge\eps_0\}.
\end{align}

\begin{lemma}\label{lem_A_S_X_close}
For all large $n$,
for every $m\in\{1,\ldots,M\}$, one has on the event
$\Om^{n,m}$,
\begin{align}\label{eq2342}
\|\tilde \bS^n\circ T^n - \upsi^{m,2}\circ\brho\|_T\le 3\delta_1/2,\qquad \|\tilde A^n - \upsi^{m,1}\|_T\le \delta_1
\end{align}
and
\begin{align}\label{eq2343}
\|\check{\bS}^n\circ T^n - \psi^{m,2}\circ\brho\|_T<2\delta_1,
\qquad
\|\check{A} - \psi^{m,1}\|_T<2\delta_1.
\end{align}
Moreover, for all large $n$ and all $t,u\in[0,\tau^n]$ such that $|t-u|<\nu_1$,
\begin{align}\label{eq2344}
|\check{X}^n(t) - \check{X}^n(u)|\le \eps_1.
\end{align}
\end{lemma}

\begin{lemma}\label{lem_taun>T}
For all large $n$, \eqref{eq2338a} holds on the event
$\cup_{m=1}^M\Om^{n,m}$.
\end{lemma}
These two lemmas are proved at the end of the section.

\skp\noi{\bf Step 2.}
As mentioned earlier, for sufficiently small $\eps_0$ one has $a^*=\Barr$.
Moreover, by Lemma \ref{lem_taun>T} and the discussion in the beginning of
step 1, no forced rejections occur on the event under consideration.
Consider the balance equation \eqref{02}
and recall that $\check Z^n$ and $\check R^n$ are
nonnegative, nondecreasing processes. Recall $\check Y^n$ defined in \eqref{16}.
Also, these processes are flat
on the set of times where $\check X^n>0$ and $\check X^n<a^*$, respectively,
where we used work conservation and the absence of forced rejections.
Recalling from Section \ref{sec_skorohod} the characterization of the Skorohod
map on an interval, it follows that, on the event under consideration,
\begin{align}\label{03}
(\check{X}^n,\check{Z}^n, \check{R}^n)(t)&=\Gam_{[0,\beta_0]} (\check{Y}^n
+ \check{A}^n - \check{\bS}^n\circ T^n).
\end{align}
Compare this relation with \eqref{eq2326}.
Let $n$ be sufficiently large so that
$\|\check{Y}^n-\by\|_T\le \delta_1$. Then by \eqref{eq2343},
\begin{align}\notag%\label{eq2376}
\left\|\left(\check{Y}^n  + \check{A}^n - \check{\bS}^n\circ T^n\right)
- \left(\by+\psi^{m,1}-\psi^{m,2}\right)\right\|_T\le 5\delta_1.
\end{align}
From the above, using \eqref{eq2316b}, \eqref{eq2326}, and Lemma \ref{lem_appendix1} it follows that on the event
$\Om^{n.m}$ one has
\begin{align}\label{eq2377}
\|\check{R}^n(T)-\vr^m(T)\|\leq \eps_0 \quad\text{and}\quad\|\check{X}^n-\ph^m\|_T
\leq \osc_D(\eps_0,h).
\end{align}
Now we bound $H^n_T$. Let $L=\sum_ih_i$.
For sufficiently large $n$,
\begin{align}\notag%\label{eq2377a}
&\int_0^T\uh\cdot\tilde X^n(u)du \le \int_0^T\uh\cdot\gamma^a( \check{X}^n(u))du+2L\eps_0 T = \int_0^Th^a( \check{X}^n(u))du+2L\eps_0 T\\\notag
&\quad\le \int_0^T h( \check{X}^n(u))du + (2L+1)\eps_0 T+T\om_1(\eps_0)
 \le \int_0^T h(\ph^m(u) )du + (2L+2)T\eps_0+T\om_1(\eps_0),
\end{align}
where the first inequality follows by Lemma \ref{lem_taun>T}, the equality follows by the definitions of $h^a$ and $\gamma^a$, the second inequality follows by \eqref{eq2307b}, and the last inequality follows by \eqref{eq2377}.
Also notice that for sufficiently large $n$
\begin{align}\notag%\label{eq2377b}
\ur\cdot\tilde R^n=r_{i^*}\tilde R^n_{i^*}=\frac{r}{\mu_{i^*}}\tilde R^n_{i^*}\le r(1+\eps_0)\theta^n_{i^*}\tilde R^n_{i^*}=r(1+\eps_0)\utheta^n\cdot\tilde R^n=r(1+\eps_0)\check{R}^n,
\end{align}
where the first and third equalities follow from \eqref{eq2340}, the second equality follows since $r=r_{i^*}\mu_{i^*}$, the inequality follows since $\utheta^n\to\utheta$, and the last equality follows by the definition of $\check{R}^n$. Denote
$\om_2(\eps_0)=((2L+2)T+r(1+\eps_0))\eps_0+T\om_1(\eps_0)$. Then from \eqref{eq2307b} one has $\om_2(0+)=0$. Moreover,
\begin{align}\label{eq2378h}
&\E[e^{b_n^2[H^n_T\wedge {J}_1]}1_{\Om^{n,m}}]\\\notag
&\leq \E[e^{b_n^2[\int_0^T\uh\cdot\tilde X^n(u)du+\ur\cdot\tilde R^n(T)]}1_{\Om^{n,m}}]\\\notag
&\le \E[e^{b_n^2[\int_0^T h(\ph^m(u) )du+r(1+\eps_0)\vr^m(T)+\om_2(\eps_0)]}1_{\Om^{n,m}}]\\\notag
&\leq\max_{0\leq t\leq T}e^{b_n^2[\int_0^t h(\ph^m(u) )du+r(1+\eps_0) \vr^m(t)+\om_2(\eps_0)]}\PP[(\tilde{A}^n,\tilde{S}^n)\in \A^m_t].
\end{align}
Let $t^m$ be such that
\begin{align}\label{eq2379}
&\max_{0\leq t\leq T}\left[e^{b_n^2[\int_0^t h(\ph^m(u) )du+r(1+\eps_0)\vr^m(t)+\om_2(\eps_0)]}\PP[(\tilde{A}^n,\tilde{S}^n)\in \A^m_t]\right]\\\notag
&\quad\leq
e^{b_n^2[\int_0^{t^m}h(\ph^m(u) )du+r(1+\eps_0)\vr^m(t^m)+\om_2(\eps_0)]}\PP[(\tilde{A}^n,\tilde{S}^n)\in \A^m_{t^m}]+\eps_0.
\end{align}
By Proposition \ref{moderate} and \eqref{eq2323},
\begin{align}\notag%\label{eq2380}
\frac{1}{b_n^2}\log\PP((\tilde{A}^n,\tilde{S}^n)
\in \oo{\A^m_{t^m}})\leq -\inf_{\psi\in\oo{\A^m_{t^m}}}\Jr(t^m,\psi)+\eps_0\le -\Jr(t^m,\upsi^m)+2\eps_0\le-\Ir(t^m,\psi^m)+2\eps_0.
\end{align}
This, along with \eqref{eq2379} show that the r.h.s.\ of \eqref{eq2378h} is
bounded by
\begin{align}\notag%\label{eq2381}
e^{b_n^2[\int_0^{t^m}h(\ph^m(u))du+r(1+\eps_0)\vr^m(t^m)-\Ir(t^m,\psi^m)+\om_2(\eps_0)]}+\eps_0.
\end{align}
Thus \eqref{eq2332} follows.
\qed

\skp

\noi{\bf Proof of Lemma \ref{lem_A_S_X_close}:}
Recall relation \eqref{eq2013} and the fact that $\tilde X^n(t)$ remains bounded.
This, along with \eqref{eq2329}
give $\|\tilde Z^n\|\leq C_{4}(1+\La^n)$.
Denote
\begin{align}\label{eq2348b}
\La_{J}=\sup_{\psi\in\AC_{J}}\|\psi^{1}\|_T + \|\psi^{2}\|_T<\iy,
\end{align}
where the finiteness follows by \eqref{eq2123} and the definition of $\AC_{J}$.
Thus, on the event $\cup_m\Om^{n,m}$, one has
$\|\tilde Z^n\|_T\leq \La_{J}+2\delta_1$.
Recalling the expression \eqref{eq2014} for $\tilde Z^n$,
it follows that $\|\brho-T^n\|_T<\nu_1$ for large $n$.
Therefore, for every $m\in\{1,\ldots,M\}$,
$\|\upsi^{m,2}\circ T^n-\upsi^{m,2}\circ\brho\|_T<\delta_1/2$.
Moreover, on the event $\Om^{n,m}$,
$\|\tilde \bS^n\circ T^n-\upsi^{m,2}\circ T^n\|_T<\delta_1$, and therefore
\[
\|\tilde \bS^n\circ T^n - \upsi^{m,2}\circ \brho\|_T
\le 3\delta_1/2.
\]
Similarly, $\|\tilde A^n - \upsi^{m,1}\|_T\le \delta_1$.

By \eqref{eq2342}, \eqref{eq2317}, and since $\utheta^n\to\utheta$, it follows that on the event $\Om^{n,m}$, \eqref{eq2343} holds.

It remains to prove \eqref{eq2344}. Fix $0\le u\le t <\tau^n$ such that $t-u<\nu_1$.
By the definition of the time $\tau^n$, an argument as that leading to \eqref{03}
shows
\begin{align}\notag%\label{eq2352}
(\check{X}^n,\check{Z}^n,\check{R}^n)(u) = \Gam_{[0,\Barr]}[\check{W}^n](u),\quad u\in[0,\tau^n).
\end{align}
where
\begin{align}\notag%\label{eq2353}
\check{W}^n(u)=\check{X}^n(0)+\check{y}^nu+\check{A}^n(u)-\check{S}^n(T^n(u))-\check{R}^n(u).
\end{align}
If we show that
$|\check{W}^n(t)-\check{W}^n(u)|<11\delta_1$ then the result follows by \eqref{eq2316b}.
Using \eqref{eq2343} along with \eqref{eq2317} and the fact
$$\limn\|\check{Y}^n - \by\|_T=0$$
shows that, for large $n$,
\begin{align}\notag%\label{eq2355}
|\check{W}^n(t)-\check{W}^n(u)|&\le  |\check{W}^n(t)-(x+yt+\psi^{m,1}(t)-\psi^{m,2}(t))|\\\notag
&\quad+|\check{W}^n(u)-(x+yu+\psi^{m,1}(u)-\psi^{m,2}(u))|\\\notag
&\quad+|\psi^{m,1}(t)-\psi^{m,1}(u))|
+|\psi^{m,2}(t)-\psi^{m,2}(u))|
< 11\delta_1.
\end{align}
This completes the proof.
\qed

\skp\noi{\bf Proof of Lemma \ref{lem_taun>T}:}
The structure of the proof borrows ideas from the proof of Lemma 4.1 of \cite{ata-shi}
(however, the content is different, as \cite{ata-shi} addresses weak convergence).
We begin with the case where the initial state lies
close to the minimizing curve. That is,
\begin{equation}\label{21}
\max_{i}|\Del^n_i(0)|\le\eps_0.
\end{equation}
At the last step of the proof we relax this assumption.

By reducing to a subsequence of $\{n\}$, one then has
that there exists an $m\in\{1,\ldots,M\}$ such that on the event
$(\tilde A^n,\tilde S^n)\in\A^m_T$ one has $\tau^n \le T$.
Let $j=j^n$ and $\xi^n$ be the corresponding
components from the representation $(j,\xi)$ of $\check{X}^n(\tau^n)$ (with $w=\check{X}^n(\tau^n)$).
Fix a positive integer $K=K(\eps_1)=[D/\eps_1]$, where
$\eps_1\le \theta_{\min}\eps_0/8$ as defined right before \eqref{eq2316b}.
Consider the covering of $[0,D]$ by the
$K-1$ intervals $\X_k=\bB(k\eps_1,\eps_1)$, $k=1,2,\ldots,K-1$, where $\bB(x,a)$ denotes
$[x-a,x+a]$. Let $\tilde\X_k=\bB(k\eps_1,2\eps_1)$.
From \eqref{eq2344} we obtain
\begin{equation}\label{eq2357}
\text{if $\check{X}^n(\tau^n)\in\X_k$ then $\check{X}^n(t)\in\tilde\X^k$ for every $t\in\bT^n:=[((\tau^n-\nu_1)\vee 0),\tau^n]$}.
\end{equation}
By considering a further subsequence, we may assume that there exists
a $k=k(m)$ such that $\check{X}^n(\tau^n)\in\X_k$ for all $n$.
The value assigned by the policy to $B^n$ (see \eqref{eq2311}) remains fixed as
$\check{X}^n$ varies within any of the intervals $(\hat a_j,\hat a_{j-1})$.
Aiming at showing a contradiction for each $k$,
we consider the following four cases.
\\ (i) $\tilde\X_k\subset(0,a^*)$ and for all $j$,
$\hat a_j\notin\tilde\X_k$.
\\ (ii) $\tilde\X_k\subset(0,a^*)$
but $\hat a_j\in\tilde\X_k$ for some $j\in\{1,2,\ldots,I-1\}$.
\\ (iii) $0\in\tilde\X_k$.
\\ (iv) $a^*\in\tilde\X_k$.

(i) $\tilde\X_k\subset(0,a^*)$ and for all $j$, $\hat a_j\notin\tilde\X_k$.
Then all points $x$ in
$\tilde\X_k$ lead to the same $j$ in the representation $(j,\xi)$ of $x$ given by
\eqref{eq2304}.
This $j=j(k)$ depends on $k$ only, and in particular does not vary with $n$.
Fix $i>j$ (except when $j=I$).
Note that $\gamma^a_i(\check{X}^n(\tau^n))=a_i$ (because $i>j$).
We show first that for sufficiently large $n$ one has
\begin{equation}\label{eq2358}
\text{for every } i> j,\quad\text{one has }
\Del^n_i(t) <\eps_0\text{ for all } t\in[0,T].
\end{equation}
This is done as follows. Assume to the contrary that $\tau^n\le T$ and that $\Del^n_i(\tau^n) \ge \eps_0$. Then, since the  jumps of $\tilde X^n_i$ are of size $(b_n\sqrt{n})^{-1}$ it follows that there must exist
$\eta^n\in[0,\tau^n]$ with the properties that
\begin{equation}\label{eq2359}
\tilde X^n_i(\eta^n)<a_i+\eps_0/2,\qquad \tilde X^n_i(t)>a_i \text{ for all } t\in[\eta^n,\tau^n].
\end{equation}
Therefore, during the time interval $[\eta^n,\tau^n]$,
$i$ is always a member of $\calH^+(\tilde X^n)$, and therefore by
\eqref{eq2311}--\eqref{eq2312}, $B^n_i(t)=\rho'_i(\tilde X^n(t))>\rho_i+C_4$, for some constant $C_4>0$.
Thus by \eqref{eq2014},
\[
\frac{d}{dt}\tilde Z^n_i\le-\frac{\mu^n_i}{b_n\sqrt n}C_4.
\]
Moreover, if we define $\hat\eta^n=\eta^n\vee(\tau^n-\nu_1)$ then by \eqref{eq2357} we get that for all
$t\in[\hat\eta^n,\tau^n]$ one has $\check{X}^n(t)\in\tilde\X_k\subset(0,a^*)$
and therefore no rejections occur.
Using these facts in \eqref{eq2014}, we have
\begin{equation}\label{eq2360}
\tilde X^n_i(\tau^n)-\tilde X^n_i(\hat\eta^n)\le
[\tilde A^n_i(\tau^n)-\tilde A^n_i(\hat\eta^n)]
-[\tilde S^n_i(T^n_i(\tau^n))-\tilde S^n_i(T^n_i(\hat\eta^n))]
-\frac{\mu^n_i}{b_n\sqrt n}C_4(\tau^n-\hat\eta^n).
\end{equation}
Fix a sequence $r_n>0$ with $r_n\to0$ and $r_n b_n\sqrt n\to\iy$.
If $\tau^n-\eta^n<r_n$ and $n$ is sufficiently large then $\hat\eta^n=\eta^n$,
thus by the definitions of $\tau^n$ and $\eta^n$, one has $\tilde X^n_i(\tau^n)-\tilde X^n_i(\hat\eta^n)\ge\eps_0/2$.
As a result,
\begin{equation}\label{eq2361}
[\tilde A^n_i(\tau^n)-\tilde A^n_i(\eta^n)]
-[\tilde S^n_i(T^n_i(\tau^n))-\tilde S^n_i(T^n_i(\eta^n))]
\ge\eps_0/2.
\end{equation}
From Lemma \ref{lem_A_S_X_close} and inequality \eqref{eq2317} it follows that the l.h.s.~of the above is smaller than $6\del_1$. Altogether, $6\del_1>\eps_0/2$, which contradicts the choice of $\del_1$ (see \eqref{eq2316b}).
If, on the other hand, $\tau^n-\eta^n\ge r_n$ then by \eqref{eq2360},
\begin{equation}\label{eq2362}
[\tilde A^n_i(\tau^n)-\tilde A^n_i(\hat\eta^n)]
-[\tilde S^n_i(T^n_i(\tau^n))-\tilde S^n_i(T^n_i(\hat\eta^n))]
\ge C_4r_nb_n\sqrt{n}.
\end{equation}
The l.h.s.~of the above is bounded from above by $2\La^n$, which is bounded by Lemma \ref{lem_A_S_X_close} and the definition of $\La_J$, see \eqref{eq2348b}. This contradicts the
fact $C_4r_nb_n\sqrt{n}\to\iy$. Therefore, \eqref{eq2358} holds.

Next, fix $i<j$ (provided $j\ne1$). Then $\gamma^a_i(\check{X}^n(\tau^n))=0$ and
whenever $\tilde X^n_i>0$, $i$ is a member of the high priority set $\calH^+(\tilde X^n)$.
This is due to the fact that there must exist $l>i$ such that $\tilde X^n_l<a_l$;
otherwise, the workload would be at least $\sum_{p=l}^{I}a_p\theta^n_p>\check{X}^n$.
Hence, the same argument yields a contradiction. Therefore
\begin{equation}\label{eq2363}
\text{for every } i<j,\quad
\Del^n_i(t)<\eps_0\text{ for all }  t\in[0,T].
\end{equation}

Consider now $j$ itself. We will show, for the case $j<I$, that
\begin{equation}\label{eq2364}
\Del^n_j(t)<\eps_0\text{ for all }  t\in[0,T].
\end{equation}
Suppose that we show
that for every $t\in[0,\tau^n]$ and every large $n$,
\begin{equation}\label{eq2365}
\text{if $\Del^n_j(t)\in(\eps_0/2,\eps_0)$, then $j\in\calH^+(\tilde X^n(t))$.}
\end{equation}
Then, assuming to the contrary that $\tau^n\le T$ and
$\tilde X^n_j(\tau^n)\ge \gamma^a_j(\check{X}^n(\tau^n))+\eps_0,$
implies that there exists
$\eta^n\in[0,\tau^n]$ with the properties that
\begin{equation}\notag%\label{eq2365c}
\tilde X^n_j(\eta^n)-\gamma^a_j(\check{X}^n(\tau^n))<3\eps_0/4,
\end{equation}
\begin{equation}\notag
\tilde X^n_j(t)-\gamma^a_j(\check{X}^n(\tau^n))>\eps_0/2 \text{ for all } t\in[\eta^n,\tau^n].
\end{equation}
From \eqref{eq2365} during the time interval $[\eta^n,\tau^n]$,
$j$ is always a member of $\calH^+(\tilde X^n_j(t))$.
Therefore, we still have inequality \eqref{eq2360} valid.
Arguing separately for the cases $\tau^n-\eta^n<r_n$ and
$\tau^n-\eta^n\ge r_n$, leads, in analogy to \eqref{eq2361} and \eqref{eq2362}, to a contradiction. There is a slight difference in the first case. If $\tau^n-\eta^n<r_n$ and $n$ is sufficiently large then $\hat\eta^n=\eta^n$,
and now $\tilde X^n_j(\tau^n)-\tilde X^n_j(\hat\eta^n)\ge\eps_0/4 + \gamma^a_j(\check{X}^n(\tau^n))-\gamma^a_j(\check{X}^n(\hat\eta^n))$. Recall that $\tau^n-\hat\eta^n<\nu_1$. Therefore, by taking $\delta_1$ to be sufficiently small, one can verify  that $|\gamma^a_j(\check{X}^n(\tau^n))-\gamma^a_j(\check{X}^n(\hat\eta^n))|$ would be significantly smaller than $\eps_0/2$. This argument follows by the continuity of $\gamma^a$ and by similar arguments to the ones we used in order to prove \eqref{eq2344}.

Now we show that \eqref{eq2365} holds (except in the case $j=I$).
Since $\utheta\cdot\gamma^a(\utheta\cdot \bar x)=\utheta\cdot \bar x$ for all $\bar x\in\calX$, $\utheta^n\to\utheta$, $\gamma^a$ is uniformly continuous, and $\calX$ is bounded, we have
\begin{equation}\label{eq2368}
  q_n:=\sup_{\bar x\in\calX}|\utheta\cdot\gamma^a(\utheta^n\cdot \bar x)-\utheta\cdot\bar x|
  \to0,\quad\text{as } n\to\iy.
\end{equation}
Note by \eqref{eq2368} that
\(
|\utheta\cdot\tilde X^n(t)-\utheta\cdot\gamma^a(\check{X}^n(t))|\le q_n\to0.
\)

Fix $t\in[0,\tau^n]$.
If $\Del^n_j(t)\ge\eps_0/2$ then
\begin{equation}\label{eq2369}
-\theta_j\eps_0/2\ge\sum_{i\ne j}\theta_i(\tilde X^n_i(t)-\gamma^a_i(\check{X}^n(t)))-\|\utheta\|q_n
\ge\sum_{i>j}\theta_i(\tilde X^n_i(t)-a_i)-\|\utheta\|q_n,
\end{equation}
where we used $\gamma^a_i(\check{X}^n(t))=0$ for $i<j$ and $\gamma^a_i(\check{X}^n(t))=a_i$ for $i>j$. These two equations hold from \eqref{eq2357}.
For all large $n$, this implies $\tilde X^n_i(t)<a_i$ for at least one $i>j$, by which
$j\in\calH^+(\tilde X^n)$.

We can now show that for every $t\in[0,T]$ one has $\max_{i\le I}|\Del^n_i(t)|<\eps_0$. Indeed, in the case $j=I$,
we have by \eqref{eq2363},
$\max_{i< I}|\Del^n_i(t)|<\eps_0$.
By \eqref{eq2368}, $|\utheta\cdot\Del^n(t)|\le q_n$.
Since $\utheta\in(0,\iy)^I$ and $q_n\to0$, we obtain
\begin{equation}\label{eq2370}
\max_{i\le I}|\Del^n_i(t)|<\eps_0.
\end{equation}
In the case $j<I$, combining \eqref{eq2358}, \eqref{eq2363}, \eqref{eq2364},
we have $\max_{i\le I}\Del^n_i(t)<\eps_0$.
Using again the fact $|\utheta\cdot\Del^n(\tau^n)|\le q_n\to0$ shows
that \eqref{eq2370} is valid in this case as well.

(ii) $\tilde\X_k\subset(0,a^*)$ but $\hat a_j\in\tilde\X_k$ for some $j\in\{1,2,\ldots,I-1\}$.
Let $(j^n(t),\xi^n(t))$ denote the representation \eqref{eq2304} for $\check{X}^n(t)$.
Note that in the
time window $\bT^n$, $j^n$ varies between two values, namely $j$ and $j+1$,
and so it is no longer true that
$\gamma^a_{j+1}(\check{X}^n(t))=a_{j+1}$ on that time interval. However, it is true that
\begin{equation}\label{eq2371}
\gamma^a_{j+1}(\check{X}^n(t))\ge
a_{j+1}-4\eps_1/\theta_{\rm min}\ge a_{j+1}-\eps_0/2,\qquad
t\in\bT^n,
\end{equation}
where the second inequality follows since $\eps_1<\theta_{\rm min}\eps_0/8$.
Indeed, we have for any $w\in\tilde\X_k$, $|w-\hat a_j|\le 4\eps_1$, since
$\hat a_j$ is also in $\tilde\X_k$. Now, if $w\ge\hat a_j$
then $\gamma^a_{j+1}(w)=a_{j+1}$. Otherwise,
\begin{equation}\notag%\label{eq2372}
w=\hat a_{j+1}+\theta_{j+1}\xi=\hat a_j-\theta_{j+1}a_{j+1}+\theta_{j+1}\xi,
\end{equation}
hence $|a_{j+1}-\xi|\le4\theta_{j+1}^{-1}\eps_1$ and \eqref{eq2371} follows.

By the same arguments as in case (i) we get that
\begin{equation}\notag
\text{for every } i\ne j+1,\quad\text{one has }
\Del^n_i(t) <\eps_0\text{ for all } t\in[0,T].
\end{equation}

As for $i=j+1$, assume to the contrry that $\tau^n\le T$ and that $\Del^n_{j+1}(\tau^n) \ge\eps_0$. Then by \eqref{eq2371} we get that $\tilde X^n_{j+1}(\tau^n)\ge a_{j+1}+\eps_0/2$ and the same arguments as in case (i) are valid.

Combining all the estimates gives,
\begin{equation}\notag%\label{eq2373}
\max_{i\le I}\Del^n_i(t)<\eps_0,\quad t\in[0,T].
\end{equation}
Along with the fact $|\theta\cdot\Del^n(\tau^n)|\le q_n$, this gives
\begin{equation}\notag%\label{eq2374}
\max_{i\le I}|\Del^n_i(t)|<\eps_0,\quad t\in[0,T].
\end{equation}

(iii) $0\in\tilde\X_k$.
This differs from case (i) in that during $\bT^n$, $\tilde X^n$ may hit
zero, and therefore $B^n$ might vanish. Note however that the way
case (i) is handled, one focuses only on time intervals where $\tilde X^n\ne0$,
and therefore the proof is valid here as well. We thus have
$\max_{i\le I}|\Del^n_i(t)|<\eps_0,\quad t\in[0,T]$.

(iv) $a^*\in\tilde\X_k$. In this case, $\utheta\cdot\tilde X^n$ may exceed
the threshold $a^*$, and rejections may occur.
The argument provided in case (i) is then slightly changed.
A negative term is added to the r.h.s.\ of \eqref{eq2360}, but the consequences
of \eqref{eq2360} remain valid with this addition.
(Note that for small $\eps_0$, $\hat a_i\ne a^*$ holds for all $i$, hence assuming
$\eps_0$ is small, we do not need to check case (ii) here.)

Having shown that $\max_{i\le I}|\Del^n_i(t)|<\eps_0,\; t\in[0,T]$ in all cases completes the proof of the lemma under \eqref{21}.

The relaxation of \eqref{21} is performed by showing that within a short time $t$,
$\max_i|\Del^n_i(t)|\le\eps_0$. This is sufficient, because on the remaining time interval
the argument provided above for the case \eqref{21} gives the result.

Fix $\delta>0$. We will show that for each $i$
and all sufficiently large $n$, there exists $t\in [0,\delta]$ such
that $|\Del^n_i(t)|\le \eps_0$. Since the proof provided above for the case \eqref{21}
treats each $i$ separately, this will assure that once $|\Del^n_i|$ is bounded by
$\eps_0$ for some $i$, it remains so for the remaining time interval.

We thus fix $i$ and prove that for
all sufficiently large $n$, there exists $t\in [0,\delta]$ such
that $|\Del^n_i(t)|\le \eps_0$.
Assume to the contrary that $|\Del^n_i|>\eps_0$ on $[0,\delta]$.
Since the jumps of $\Del^n$ are of order $n^{-1/2}$, we either have $\Del^n_i>\eps_0$ on $[0,\del]$,
or $\Del^n_i<-\eps_0$ on $[0,\del]$, provided $n$ is large.
In the former case, $i$ is always a member of $\calH^+(\tilde X^n(t))$, for every $t$,
and therefore by
\eqref{eq2311}--\eqref{eq2312}, $B^n_i(t)=\rho'_i(\tilde X^n(t))>\rho_i+C_4$, for some constant $C_4>0$.
Thus by \eqref{eq2014},
\[
\frac{d}{dt}\tilde Z^n_i\le-\frac{\mu^n_i}{b_n\sqrt n}C_4.
\]
Therefore, as was argued in \eqref{eq2360},  by \eqref{eq2014} we obtain that
\begin{align*}
\frac{\mu^n_i}{b_n\sqrt n}C_4\delta &\leq -[\tilde X^n_i(\delta)-\tilde X^n_i(0)]+
[\tilde A^n_i(\delta)-\tilde A^n_i(0)]
-[\tilde S^n_i(T^n_i(\delta))-\tilde S^n_i(T^n_i(0))]\\\notag
&\qquad - [\tilde R^n_i(\delta)-\tilde R^n_i(0)]
\\
& \le 2D_i+2\La^n.
\end{align*}
The r.h.s.~of the above bounded from above, as follows from
Lemma \ref{lem_A_S_X_close}
and the fact that $\La_J<\iy$ (see \eqref{eq2348b}). This contradicts the
fact that the l.h.s.~goes to infinity.

In case that $\Del^n_i<-\eps_0$ on $[0,\del]$,
$i$ is not a member of $\calH^+(\tilde X^n(t))$ for any $t$, and therefore by
\eqref{eq2311}--\eqref{eq2312}, $B^n_i(t)=0$.
Thus by \eqref{eq2014},
\[
\frac{d}{dt}\tilde Z^n_i=\frac{\mu^n_i}{b_n\sqrt n}\rho_i.
\]
Hence
\begin{align*}
\frac{\mu^n_i}{b_n\sqrt n}\rho_i\delta=& [\tilde X^n_i(\delta)-\tilde X^n_i(0)]-
[\tilde A^n_i(\delta)-\tilde A^n_i(0)]
+[\tilde S^n_i(T^n_i(\delta))-\tilde S^n_i(T^n_i(0))]\\\notag
&+ [\tilde R^n_i(\delta)-\tilde R^n_i(0)].
\end{align*}
The r.h.s.~of the above is bounded, using similar considerations along with \eqref{eq2329},
whereas again, the l.h.s.\ tends to infinity with $n$.
\qed

\appendix
\section{Appendix}\beginsec
\manualnames{A}

\subsection{The Skorohod map on a finite interval}\label{sec_skorohod}

Given $a$ and $b$, $a<b$, the {\it Skohorod map} on the interval $[a,b]$
maps $\calD([0,\iy),\R)$ to $\calD([0,\iy),\R^3)$.
It is denoted by $\Gam_{[a,b]}$, and is characterized as the solution
map $\om\to(\ph,\eta_1,\eta_2)$ to the
problem of finding, for a given $\om$, a triplet $(\ph,\eta_1,\eta_2)$, such that
\begin{equation}\notag%\label{eq903}
\ph=\om+\eta_1-\eta_2,\qquad \ph(t)\in[a,b] \text{ for all } t,
\end{equation}
$\eta_i$ are nonnegative and nondecreasing, $\eta_i(0-)=0$, and
\begin{equation}\notag%\label{eq904}
\int_{[0,\iy)}1_{(a,b]}(\ph)d\eta_1=\int_{[0,\iy)}1_{[a,b)}(\ph)d\eta_2=0.
\end{equation}
By writing $\eta_i(0-)=0$ we adopt the convention that $\eta_i(0)>0$ is regarded a jump
at zero. This convention, in conjunction with $\int_{[0,\iy)}1_{(a,b]}(\ph)d\eta_1=0$
(resp., $\int_{[0,\iy)}1_{[a,b)}(\ph)d\eta_2=0$),
means that if $\om(0)<a$ (resp., $\om(0)>b$) then $\ph(0)=a$ (resp., $b$).
If, however,
$\om(0)\in[a,b]$ then $\ph(0)=\om(0)$, and $\eta_i$ have no jump at zero.
See \cite{KLRS} for existence and uniqueness of solutions,
and continuity and further properties of the map.
In particular, we have the following.
\begin{lemma}\label{lem_appendix1}
Fix $b>0$. Then there exists a constant $\ee$ such that
for every $T>0$, $\del>0$ and $\om,\tilde\om\in\calD([0, \iy), \R)$,
\begin{equation}\notag%\label{eq2400}
\|\Gam_{[0,b]}(\om) - \Gam_{[0,b]}(\tilde\om)\|_{T}
\le \ee \|\om-\tilde\om\|_T,
\end{equation}
and
\[
\osc_T(\del,\Gam_{[0,b]}(\om))\le \ee\,\osc(\del,\om).
\]
\end{lemma}

\subsection{On the rate functions $\Ir$ and $\Jr$}

For every $T\in\R_+$ and $\psi=(\psi^1,\psi^2)\in\calP^2$ set $\tilde\Ir(T,\psi)=\tilde\Ir_1(T,\psi^1)+\tilde\Ir_2(T,\psi^2)$,
\begin{equation}\label{eq2026d}
\tilde\Ir_1(T,\psi^1) = \inf \{\Jr_1(T,\upsi^1):\upsi^1\in\calP^I, \utheta\cdot\upsi^1=\psi^1\}
\end{equation}
and
\begin{equation}\label{eq2026dd}
\tilde\Ir_2(T,\psi^2) = \inf \{\Jr_2(T,\upsi^2):\upsi^2\in\calP^I, \utheta\cdot(\upsi^2\circ\brho)=\psi^2\},
\end{equation}

Recall that $\Ir$ is defined in \eqref{eq2026f} and $\Jr$ in \eqref{eq2018}.
\begin{lemma}\label{lem_underline_psi}
For every $\psi=(\psi^1,\psi^2)\in\calP^2$ there exists
$\upsi=(\upsi^1,\upsi^2)\in\calP^{I}\times\calP^{I}$ such that
\begin{equation}\label{eq2026e}
(\utheta\cdot\upsi^1,\utheta\cdot(\upsi^2\circ\brho))=(\psi^1,\psi^2)
\end{equation}
and
\begin{equation}\label{eq2026ee}
\tilde\Ir_k(T,\psi^k)=\Ir_k(T,\psi^k) = \Jr_k(T,\upsi^k),\quad k=1,2, \,T>0.
\end{equation}
\end{lemma}

\noi
{\bf Proof:}
Define $\Kr:[0,\iy)\times\calP^I\to[0,\iy]$ by
\begin{align}\notag%\label{app00}
\Kr(T,\upsi)=\left\{\begin{array}{ll}
               \sum_{i=1}^I \al_i
               \int_0^T\dot\psi_i^2(u)du &\ \mbox{if all}\ \psi_i\in\AC_0([0,T],\R),
\\
                \infty & \ \mbox{otherwise},
              \end{array}
\right.
\end{align}
where $\al_1,\ldots,\al_I>0$. Then both $\Jr_1$ and $\Jr_2$ are of the form $\Kr$.
Set $l_1,\ldots,l_I\in(0,1]$. Define $\Lr:[0,\iy)\times\calP\to[0,\iy]$ by
\begin{equation}\label{app01}
\Lr(T,\psi) = \inf \{\Kr(T,\upsi):\upsi\in\calP^I, \utheta\cdot\upsi(l)=\psi\},
\end{equation}
where $\upsi(l)=(\psi_1(l_1t),\ldots,\psi_I(l_It))$. Note that in case $l_1=\ldots=l_I=1$ (resp., $l_1+\ldots+l_I=1$) the function $\Lr$ gives $\tilde\Ir_1$ (resp., $\tilde\Ir_2$).
We show that for every $\psi\in\calP$ there is $\upsi$ such that
\begin{equation}\label{app02}
\utheta\cdot\upsi(l)=\psi
\end{equation}
and
\begin{equation}\label{app03}
\Lr(T,\psi) = \Kr(T,\upsi).
\end{equation}
We now calculate $\Lr$:
\begin{align}\notag%\label{app04}
&\inf \left\{\Kr(T,\upsi):\upsi\in\calP^I, \utheta\cdot\upsi(l)=\psi\right\}\\\notag
&\quad=
\inf\left\{\sum_{i=1}^I\al_i\int_0^T\dot\psi_i^2(u)du :
\sum_{i=1}^I\theta_i\psi_i(l_iu)=\psi(u),\;u\in[0,T]\right\}\\\notag
&\quad=
\inf\left\{\sum_{i=1}^I\al_i\int_0^{l_iT}\dot\psi_i^2(u)du :
\sum_{i=1}^I\theta_i\psi_i(l_iu) =\psi(u),\;u\in[0,T]\right\}\\\notag
&\quad=
\inf\left\{\sum_{i=1}^I\al_il_i\int_0^T\dot\psi_i^2(l_it)dt:
\sum_{i=1}^I\theta_i\psi_i(l_iu) =\psi(u),\;u\in[0,T]\right\}\\\notag
&\quad=
\inf\left\{\sum_{i=1}^I\frac{\al_i}{l_i}\int_0^T\left(\frac{d}{dt}\psi_i(l_it)\right)^2
dt : \sum_{i=1}^I\theta_i\psi_i(l_iu) =\psi(u),\;u\in[0,T]\right\}.
\end{align}
Above, the second equality follows since for every $i\in\calI$ the constraint on $\psi_i$ exists only on the time interval $[0,l_iT]$. Hence, we are free to choose $\psi_i$ on the interval $(l_iT,T]$. By simple calculation, the infimum on the r.h.s.\ is attained by $\upsi=(\psi_1,\ldots,\psi_I)$,
\begin{align}\notag%\label{app06}
\psi_i(u)=\left\{\begin{array}{ll}
               \frac{\theta_il_i}{\al_i}\left(\frac{\theta_1^2l_1}{\al_1}+\cdots+\frac{\theta_I^2l_I}{\al_I}\right)^{-1}\psi(u/l_i) &\ 0\le u\le l_iT,
\\
                \psi_i(l_iT) & \ u>l_iT.
              \end{array}
\right.
\end{align}
By substitution,
\begin{align}\notag%\label{app07}
\Kr(T,\upsi)=\left(\frac{\theta_1^2l_1}{\al_1}+\cdots+\frac{\theta_I^2l_I}{\al_I}\right)^{-1}\int_0^T\dot\psi^2(t)dt.
\end{align}
By the taking proper $l_i$'s and $\al_i$'s for each of the cases $\Ir_1$ and $\Ir_2$,
one obtains \eqref{eq2026ff}. This completes the proof.
\qed

\skp

\noi{\bf
Acknowledgement.}
The authors thank an associate editor and a referee
for their helpful comments that considerably helped improve the exposition.

\footnotesize

\bibliographystyle{plain}

\bibliography{refs}
%\bibliography{Refs_Rami}

\end{document}